\documentstyle{amsppt}
\magnification=1095

\NoBlackBoxes

\topmatter
\nopagenumbers
\font\titlefont=cmbx12\font\authorfont=cmr10
\baselineskip=13pt
\rightheadtext{HERMITIAN CHARACTERISTICS OF NILPOTENT ELEMENTS}
\leftheadtext{E. TEVELEV}
\topmatter
\centerline{\titlefont
HERMITIAN CHARACTERISTICS}
\centerline{\titlefont
OF NILPOTENT ELEMENTS}
\vskip 1.2pc
\centerline{\authorfont
E. TEVELEV
}
\vskip 1.1pc
\centerline{University of Texas in Austin}
\smallskip
\centerline{tevelev\@math.utexas.edu}

\abstract
We define and study several equivariant stratifications
of the isotropy and coisotropy representations
of a parabolic subgroup in a complex reductive group.
\endabstract

\endtopmatter
\let\phi=\varphi

\def\cal#1{{\fam2 #1}}
\def\C{{\Bbb C}}

\def\R{{\Bbb R}}

\def\ov{\overline}

\def\g{\goth{g}}
\def\h{\goth{h}}
\def\z{\goth{z}}
\def\ssl{\goth{sl}}

\def\gl{\goth{gl}}

\def\k{\goth{k}}
\def\m{\goth{m}}
\def\p{\goth{p}}

\def\l{\goth{l}}
\def\n{\goth{n}}

\def\P{\Bbb P}
\def\p{\goth{p}}
\def\u{\goth{u}}
\def\End{{\text{\rm End}}}
\def\Ker{{\text{\rm Ker}}}

\def\Im{{\text{\rm Im}}}
\def\Cl{{\text{\rm Cl}}}

\def\Pf{{\text{\rm Pf}}}

\def\GL{\text{GL}}

\def\SL{\text{SL}}

\def\SO{\text{SO}}
\def\OOO{\text{O}}

\def\Sp{\text{Sp}}
\def\Spin{\text{Spin}}
\def\ad{\text{ad}}
\def\Ad{\text{Ad}}

\def\Z{\Bbb Z}
\def\O{\cal O}

\def\rk{\text{rk}}
\def\Hom{\text{Hom}}
\def\Id{\text{Id}}
\def\ZZ{\widehat Z}
\def\W{\mathop{\circ}}
\def\B{\mathop{\bullet}}
\def\Es#1#2#3#4#5#6#7{\displaystyle
#1\!\!-\!\!#2
\!\!-\!\!#3 
\!\!-\!\!#4\!\!-\!\!#5\!\!-\!\!#6
\hskip-24.2pt\lower4.5pt\hbox{${\scriptstyle|}
\hskip-3.35pt\lower6pt\hbox{$#7$
\lower3pt\hbox{\ }
}$}}
\def\Ee#1#2#3#4#5#6#7#8{\displaystyle
#1\!\!-\!\!#2
\!\!-\!\!#3 
\!\!-\!\!#4\!\!-\!\!#5\!\!-\!\!#6\!\!-\!\!#7
\hskip-24.2pt\lower4.5pt\hbox{${\scriptstyle|}
\hskip-3.35pt\lower6pt\hbox{$#8$
\lower3pt\hbox{\ }
}$}}
\def\Ess#1#2#3#4#5#6#7#8{\displaystyle
{#1}\!\!-\!\!{#2}
\!\!-\!\!{#3} 
\!\!-\!\!{#4}\!\!-\!\!{#5}\!\!-\!\!{#6}
\hskip-24.2pt\lower4.5pt\hbox{${\scriptstyle|}
\hskip-3.35pt\lower6pt\hbox{$
{\displaystyle{#7}}{\scriptstyle{#8}}$}$}}
\def\Ees#1#2#3#4#5#6#7#8#9{\displaystyle
{#1}\!\!-\!\!{#2}
\!\!-\!\!{#3} 
\!\!-\!\!{#4}\!\!-\!\!{#5}\!\!-\!\!{#6}\!\!-\!\!{#7}
\hskip-24.2pt\lower4.5pt\hbox{${\scriptstyle|}
\hskip-3.35pt\lower6pt\hbox{$
{\displaystyle{#8}}{\scriptstyle{#9}}$}$}}

\def\ra#1{\mathop{\rightarrow}\limits^{\scriptstyle#1}}
\def\da#1{\downarrow\hskip-3pt\lower-1pt\hbox{$\scriptstyle#1$}}
\def\sea#1{\searrow\hskip-7pt\lower-3pt\hbox{$\scriptstyle#1$}}
\def\w#1{\hbox to 9pt{#1\hfil}}
\def\LRA#1#2{\mathop{\Longrightarrow}\limits^{#1}_{#2}}

\document
\head \S0. Introduction\endhead

Let $G$ be a connected complex reductive algebraic group, $P\subset G$
a proper parabolic subgroup, $L\subset P$ a Levi subgroup,
$Z\subset L$ the connected component of the center of $L$,
$P^-$ the opposite parabolic subgroup, so $P\cap P^-=L$.
Let $\g$, $\p$, $\l$, $\z$, and $\p^-$ denote the corresponding Lie algebras.
Then $\g$ is $\ZZ$-graded, where $\ZZ$ is the character group of $Z$ written 
additively, 
$\g=\mathop\oplus_{\chi\in\ZZ}\g_\chi$, $\g_0=\l$.
$P$ is maximal if and only if $\ZZ$-grading
reduces to $\Z$-grading of
the form $\mathop{\oplus}\limits_{i=-k}^k\g_i$. 
Let $\n\subset\p$ be the unipotent radical of $\p$.
The group $\ZZ$ admits a total ordering such that 
$\p=\mathop{\oplus}_{\chi\ge0}\g_\chi$,
$\n=\mathop\oplus_{\chi>0}\g_\chi$.
We shall also use the identification
$\g/\p^-=\mathop\oplus_{\chi>0}\g_\chi$,
which is an isomorphism of $L$-modules.
Any element $e\in\n$ (resp.~$e\in\g/\p^-$)
has the weight decomposition $e=\sum_{\chi>0}e_\chi$, where $e_\chi\in\g_\chi$.
Each $e_\chi$ is a homogeneous nilpotent element of $\g$,
and therefore  can be embedded 
(not uniquely) in a homogeneous $\ssl_2$-triple
$\langle e_\chi,h_\chi,f_\chi\rangle$, where 
$f_\chi\in\g_{-\chi}$, $h_\chi\in\l$.
Here by an 
$\ssl_2$-triple 
$\langle e,f,h\rangle$ we mean a collection of possibly zero vectors
such that
$[e,f]=h$, $[h,e]=2e$, $[h,f]=-2f$.
The collection of elements
$\{h_\chi\}$ is called {\it a multiple characteristic\/} of~$e$
(or of the collection $\{e_\chi\}$).
A compact real form of the Lie algebra
of a reductive subgroup in $G$ is, by definition,
the Lie algebra of a compact real form of
this complex algebraic subgroup.  
We fix a compact real form $\k\subset\l$.
The multiple characteristic $\{h_\chi\}$ is called {\it Hermitian\/}
if any $h_\chi\in i\k$.

In this paper we consider the following question:
is it true that any $e\in\n$ {\rm(}resp.~$e\in\g/\p^-${\rm)}
admits a Hermitian multiple characteristic up to $P$-conjugacy
(resp.~up to $P^-$-conjugacy)? 
Clearly,  we may (and shall) suppose that 
$G$ is simple. 
We can identify parabolic subgroups with coloured Dynkin
diagrams. Black vertices correspond to simple roots such that
the corresponding root subspaces belong to the Levi part of the 
parabolic subgroup.
The following theorem is the main result of this paper.

\proclaim{Theorem 1}
Suppose that $G$ is simple and not $E_7$ or $E_8$,
or $G$ is equal to $E_7$ or $E_8$, but $P$ is not of type
38--59 {\rm(}see the Table at the end of this paper{\rm)}.
Then for any $e\in\n$ {\rm(}resp.~$e\in\g/\p^-${\rm)}
there exists $p\in P$ {\rm(}resp.~$p\in P^-${\rm)}
such that $\Ad(p)e$ admits a Hermitian multiple characteristic. 
\endproclaim

The remaining cases 38--59 will be studied in the sequel of this paper.
We remark here that the {\it adjoint case}
$e\in\g/\p^-$ and the {\it coadjoint case} $e\in\n$
are indeed being considered separately. 

Theorem~1  
has the following corollary valid over arbitrary algebraically
closed field $k$ of characterictic $0$.
For any collection of elements $v_1,\ldots,v_n$ of an algebraic
Lie algebra $\h$ we denote by $\langle v_1,\ldots,v_n\rangle_{alg}$
the minimal algebraic Lie subalgebra of $\h$ containing $v_1,\ldots,v_n$.

\proclaim{Theorem 2}
For the pairs $(G,P)$ satisfying Theorem~1 over $\C$,
the following is true over $k$.
For any $e\in\n$ {\rm(}resp.~$e\in\g/\p^-${\rm)}
there exists $p\in P$ {\rm(}resp.~$p\in P^-${\rm)}
such that $\Ad(p)e$ admits a multiple characteristic
$\{h_\chi\}$ with reductive 
$\langle h_\chi\rangle_{alg}$.
\endproclaim

The following Theorem~3 (conjectured in \cite{Te2})
was already proved in \cite{Te1}
for all parabolic subgroups that satisfy Theorem~2.
Let $\goth R_G$ denote the set of irreducible representations of $G$.
For any $V\in\goth R_G$, there exists a unique maximal proper $P$-submodule
$M_V\subset V$. Therefore, we have a linear map
$\Psi_V:\,\g/\p\to\Hom(M_V,V/M_V)$,
$\Psi_V(x)v=x\cdot v\mod M_V$.

\proclaim{Theorem 3}
For the pairs $(G,P)$ satisfying Theorem~2,
there exists an algebraic $P$-invariant stratification
$\g/\p=\mathop\sqcup_{i=1}^NX_i$
such that for any $V\in\goth R_G$,
the function $\rk\,\Psi_V$ is costant along each~$X_i$.
In other words, the linear span of functions $\rk\,\Psi_V$
in the algebra of constructible functions on $\g/\p$
is finite-dimensional.
\endproclaim

If the stratification is known explicitly
it allows to solve the following classical geometric problem effectively.
For any irreducible equivariant spanned vector bundle $\cal L$ on $G/P$,
to check whether its generic global section has zeros.
This class of geometric problems includes, for example,
the exact estimates on the maximal possible dimension
of a projective subspace in a generic projective hypersurface
of given degree, or of
an isotropic subspace of a generic skew-symmetric
form of given degree.
The corresponding algorithm and examples were given in \cite{Te2}.

The stratification from Theorem~3 provides an alternative
to the orbit decomposition for the action of $P^-$ on $\g/\p^-$.
It is known that this action has an open orbit, see \cite{LW},
however the number of orbits is usually infinite,
see \cite{PR}, \cite{BHR}.

All necessary facts about complex and real Lie groups, 
Lie algebras, and algebraic groups used in this paper without
specific references can be found in~\cite{VO}.
In particular, we numerate simple roots of simple Lie algebras
as in \cite{loc. cit.}.
If $\g$ is a simple group of rank $r$ then $\alpha_1,\ldots,\alpha_r$
denote its simple roots, $\phi_1,\ldots,\phi_r$ denote its fundamental weights.
For any dominant weight $\lambda$ we denote by $R(\lambda)$
the irreducible representation of highest weight $\lambda$.

This paper was written during my stay at
the Erwin Shr\"odinger Institute in Vienna and 
at the University of Glasgow.
I would like to thank my hosts for the warm hospitality. 

\head \S1. Hermitian Characteristics in Classical Groups\endhead

In this section we 
give some background information, 
prove Theorem~1 for classical groups, and prove Theorem~2.
Let us consider the case 
$\g=\gl(V)$ first.
We recall the notion of the Moore--Penrose inverse
of linear maps. 
Let $V_1$ and $V_2$ be vector spaces with
Hermitian scalar products.
For any linear map $F:\,V_1\to V_2$ its Moore-Penrose
inverse is a linear map $F^+:\,V_2\to V_1$ defined as follows.
Consider $\Ker\, F\subset V_1$ and $\Im F\,\subset V_2$.
Let  $\Ker^\perp F\subset V_1$ and $\Im^\perp F\subset V_2$
be their orthogonal complements with respect 
to the Hermitian scalar products.
Then $F$ defines via restriction a bijective linear map
$\tilde F:\,\Ker^\perp F\to\Im\, F$.
Then $F^+:\,V_2\to V_1$ is a unique linear map such that
$F^+|_{\Im^\perp F}=0$ and 
$F^+|_{\Im F}=\tilde F^{-1}$.
It is easy to see that 
$F^+$ is the unique solution of the following system
of {\it Penrose equations}:
$$F^+FF^+=F^+,\quad FF^+F=F,\quad FF^+\ \text{and}\ F^+F
\ \text{are Hermitian}.\eqno(*)$$
For example, if $F$ is invertible then $F^+=F^{-1}$.

Suppose now that $\g=\gl(V)$ and take
any decomposition $V=V_1\oplus\ldots\oplus V_k$.
Then we have a block decomposition 
$\gl(V)=\mathop\oplus_{i,j=1}^k\Hom(V_i,V_j)$.
The parabolic subalgebra
$\p=\mathop\oplus_{i\le j}\Hom(V_i,V_j)$ has
the unipotent radical $\n=\mathop\oplus_{i<j}\Hom(V_i,V_j)$
and any parabolic subalgebra in $\gl(V)$ has this form.
Consider the corresponding $\ZZ$-grading.
Then $\l=\gl(V)_0=\mathop\oplus_{i=1}^k\End(V_i)$ 
and all other rectangular blocks $\Hom(V_i,V_j)$ get different grades.

We claim that for any 
collection $\{e_{ij}\}$, $e_{ij}\in\Hom(V_i,V_j)$, $i<j$,
there exists a unique Hermitian multiple characteristic
with respect to any compact 
real form of~$\l$.
Indeed, we fix Hermitian
scalar products in $V_1,\ldots,V_k$
and set 
$\k=\mathop\oplus_{i=1}^k\u(V_i)$,
where $\u(V_i)$ is the Lie algebra of skew-Hermitian operators.
Any compact real form of $\l$ has this presentation.
Now let us take Moore--Penrose inverses 
$f_{ji}=e_{ij}^+\in\Hom(V_j,V_i)$
of elements $e_{ij}$.
The equations $(*)$ are equivalent to the statement
that $\langle e_{ij}, h_{ij}=[e_{ij},f_{ji}], f_{ji}\rangle$
is an $\ssl_2$-triple containing $e_{ij}$
with a Hermitian characteristic~$h_{ij}$.
Therefore, $\{h_{ij}\}_{i<j}$ is a unique Hermitian multiple characteristic
of $\{e_{ij}\}$. In particular, Theorem~1 is proved for $\gl(V)$.
The same argument applies for $\ssl(V)$.

\bigskip

\definition{Definitions}
\item{$\bullet$}
An $\Ad(L)$-orbit $\O\subset\g_\chi$ is called
an {\it ample orbit\/}
if for any $e_\chi\in\O$
there exists
a homogeneous $\ssl_2$-triple $\langle e_\chi,h_\chi,f_\chi\rangle$ 
with Hermitian characteristic $h_\chi\in i\k$.
\item{$\bullet$}
An element $x\in\g_{\chi}$ is called ample if its $L$-orbit is ample.
\item{$\bullet$}
The $\ZZ$-grading of $\g$ is called {\it ample in degree $\chi$}, if
all $\Ad(L)$-orbits in $\g_\chi$ are ample.
\item{$\bullet$}
For some parabolic subgroups the element $p$
in Theorem~1 can be chosen within the Levi subgroup $L$.
These parabolic subgroups are called {\it weakly ample}.
\enddefinition

A non-zero weight $\chi\in\ZZ$ is called {\it reduced} if $\g_{2\chi}=0$.

\proclaim{Basic Lemma {\rm(see \cite{Te1})}}
If $\chi\in\ZZ$ is reduced then the grading is ample 
in degree~$\chi$.
\endproclaim

Basic Lemma shows that almost all components of the grading are
automatically ample.
Now we can prove Theorem 1 for orthogonal and symplectic groups.

\proclaim{Theorem 4}
Any parabolic subgroup in $\SO_n(\C)$ or $\Sp_n(\C)$ is weakly
ample.
\endproclaim

Suppose that $V=\C^n$ is a complex vector space endowed with a non-degenerate
bilinear form $\omega$, which is either symmetric or skew-symmetric.
We denote by $G(\omega)\subset\SL(V)$
the corresponding special orthogonal or symplectic group
of automorphisms of $V$ preserving $\omega$. Let
$\g\subset\ssl(V)$ be its Lie algebra
of skew-symmetric operators with respect to $\omega$.
The subspace $U\subset V$ is called
isotropic if $\omega$ vanishes on $U$.
Let
$0=F_0\subset F_1\subset\ldots\subset F_k\subset V$
be a flag of isotropic subspaces. Then the stabilizer $P$
of this flag in $G(\omega)$ is a parabolic subgroup
and any parabolic subgroup has this form.

In order to fix the Levi subgroup of $P$, we choose
subspaces $U_k^+\subset F_k$ complement to $F_{k-1}$ and
we choose an isotropic subspace $G_k$ transversal to $F_k$ such that
the restriction of $\omega$ on $F_k\oplus G_k$ is non-degenerate.
Then the subgroup $L$ of $G(\omega)$ preserving
both $U_1^+,\ldots,U_k^+$ and $G_k$ is a Levi subgroup of $P$.
$F_k$ and $G_k$ are naturally dual to each other
with respect to the bilinear form:
any $v\in G_k$ defines a linear functional
$\omega(v,\cdot)$ on $F_k$.
Let
$G_k=U_1^-\oplus\ldots\oplus U_k^-$
be the decomposition dual to the decomposition
$F_k=U_1^+\oplus\ldots\oplus U_k^+$. Let $W=(F_k\oplus G_k)^\perp$.
Then $L$ automatically preserves $U_1^-,\ldots,U_k^-$ and $W$.
The semisimple part of $\l$ is isomorphic to
$\ssl(U_1^-)\oplus\ldots\oplus\ssl(U_k^-)\oplus\g(W)$, where
$\g(W)$ is the Lie algebra of skew-symmetric operators of $W$.
Each $\ssl(U_i^-)$ acts naturally in $U_i^-$, dually on $U_i^+$, and
trivially on $W$ and on $U_j^+$, $U_j^-$ for $j\ne i$.
$\g(W)$ acts naturally on $W$ and trivially on $F_k\oplus G_k$.
The center $\z$ of $\l$ acts on $U_k^+$ by scalar transformations $\lambda_k E$,
on $U_k^-$ by $-\lambda_k E$, and trivially on $W$.
Here $\lambda_1,\ldots,\lambda_k$ is a basis of $\z^*$.

The choice of a compact form in~$\l$ is equivalent to the choice
of a compact form in each $\ssl(U_i^-)$ (i.e.~the choice
of a Hermitian scalar product in $U_i^-$) and
a compact form $\k(W)$ of $\g(W)$.
We shall vary Hermitian scalar products in $U_i^-$ later,
but the compact real form of $\g(W)$ is to be fixed forever now.
We fix a basis $b_1,\ldots,b_m$ of $W$ such that the matrix of $\omega$
in this basis is equal to~$I$, where $I=\Id$ in the symmetric case
and $I=\left(\matrix0&\Id\cr -\Id&0\cr\endmatrix\right)$
in the skew-symmetric case.
We fix a standard Hermitian form in~$W$ 
by formula $\{u,v\}=\omega(u,I^t\ov v)$. 
Then the subalgebra $\k(W)$ of skew-Hermitian operators in $\g(W)$
is its compact real form.

Non-trivial $\ZZ$-graded components of $\g$ can be described as follows.

\item{$\bullet$}
Any linear operator $A_{ij}:\,U_i^+\to U_j^+$, $i\ne j$ 
gives by duality
the linear operator $A_{ij}':\, U_j^-\to U_i^-$. We define the skew-symmetric
operator $\tilde A_{ij}$ by setting
$\tilde A_{ij}|_{U_i^+}=A_{ij}$, \quad
$\tilde A_{ij}|_{U_j^-}=-A_{ij}'$,
the restriction of $\tilde A_{ij}$ on other components of the decomposition
$V=\mathop\oplus\limits_{k}U_k^+
\oplus\mathop\oplus\limits_{k}U_k^-\oplus W$ is trivial.
These operators form $\g_{\lambda_j-\lambda_i}$, which
belongs to $\n$ if $i>j$.

\item{$\bullet$}
Any linear operator $A_i':\,W\to U_i^+$ (resp.~$A_i':\,W\to U_i^-$)
determines by duality
the linear operator $A_i:\, U_i^-\to W$ (resp.~$A_i:\, U_i^+\to W$), where
we identify $W$ and $W^*$, $w\mapsto\omega(w,\cdot)$. 
We define the skew-symmetric
operator $\tilde A_i$ by setting
$\tilde A_i|_{W}=A'_i$, $\tilde A_i|_{U_i^\mp}=-A_i$,
the restriction of $\tilde A_i$ on other components of the decomposition
$V=\mathop\oplus\limits_{k}U_k^+
\oplus\mathop\oplus\limits_{k}U_k^-\oplus W$ is trivial.
These operators form $\g_{\lambda_i}\subset\n$
(resp.~$\g_{-\lambda_i}$).

\item{$\bullet$}
Any linear operator $B_{ij}:\,U_i^-\to U_j^+$
(resp.~$B_{ij}:\,U_i^+\to U_j^-$), $i\ne j$ determines 
by duality the linear operator $B_{ij}':\, U_j^-\to U_i^+$
(resp.~$B_{ij}':\, U_j^+\to U_i^-$). 
We define the skew-symmetric
operator $\tilde B_{ij}$ with respect to $\omega$ by setting
$\tilde B_{ij}|_{U_i^\mp}=B_{ij}$, $\tilde B_{ij}|_{U_j^\pm}=-B_{ij}'$,
the restriction of $\tilde B_{ij}$ on other components of the decomposition
$V=\mathop\oplus\limits_{k}U_k^+
\oplus\mathop\oplus\limits_{k}U_k^-\oplus W$ is trivial.
These operators form $\g_{\lambda_i+\lambda_j}\subset\n$
(resp.~$\g_{-\lambda_i-\lambda_j}$).

\item{$\bullet$}
Finally, 
any skew-symmetric linear operator $B_i:\,U_i^-\to U_i^+$
(resp.~$B_i:\,U_i^+\to U_i^-$) 
defines the skew-symmetric
operator $\tilde B_i$ with respect to $\omega$:
$\tilde B_i|_{U_i^\mp}=B_i$,
the restriction of $\tilde B_i$ on other components of the decomposition
$V=\mathop\oplus\limits_{k}U_k^+
\oplus\mathop\oplus\limits_{k}U_k^-\oplus W$ is trivial.
These operators form $\g_{2\lambda_i}\subset\n$
(resp.~$\g_{-2\lambda_i}$).
This component is trivial if $\dim U_i^+=1$.

\medskip

By Basic Lemma, if  $\g_{2\chi}=0$ then for any
$e\in\g_\chi$ and for any compact form $\k\subset\l$ there
exists a homogeneous $\ssl_2$ triple $\langle e,h,f\rangle$
with a Hermitian characteristic $h\in i\k$.
Therefore, in our case it remains to prove that
for any set of elements $e_i\in\g_{\lambda_i}$, $i=1\ldots k$
there exists a compact form $\k\subset \l$ and a Hermitian multiple
characteristic $h_i\in i\k$, $i=1,\ldots,k$.

Any $e_i$ is defined by a linear map
$A_i:\,U_i^-\to W$,
$f_i\in\g_{-\lambda_i}$ is defined by a linear map 
$B_i:\,W\to U_i^-$.
Using the description of weighted components of $\g$ given above,
it is easy to see that $e_i$ and $f_i$ can be embedded in a
homogeneous $\ssl_2$-triple if and only if $A_i$ and $B_i$
satisfy the following system of matrix equations
$$2A_i=2A_iB_iA_i-(A_iB_i)^\#A_i,\quad 2B_i=2B_iA_iB_i-B_i(A_iB_i)^\#,$$
where for any $A\in\Hom(W,W)$ we denote by $A^\#$
its adjoint operator with respect to $\omega$.
Moreover, the characteristic of this $\ssl_2$-triple will be Hermitian 
if and only if
$$B_iA_i,\  A_iB_i-(A_iB_i)^\#\ \text{are Hermitian  operators}.$$
 
Therefore, it remains to prove the following lemma.

\proclaim{Lemma}
Suppose that $U=\C^n$, $W=\C^k$ are complex vector spaces.
Let $\omega$ be a non-degenerate symmetric {\rm(}resp.~skew-symmetric{\rm)} form
on $W$ with matrix $I=\Id$, resp.~$I=\left(\matrix0&\Id\cr -\Id&0\cr\endmatrix\right)$.
We fix a standard Hermitian form in~$W$. 
Let $A\in\Hom(U,W)$.
Then there exists a Hermitian form on $U$ and
an operator $B\in\Hom(W,U)$ such that
$$2A=2ABA-(AB)^\#A,\quad 2B=2BAB-B(AB)^\#,\eqno(*)$$
$$BA,\  AB-(AB)^\#\ \text{are Hermitian  operators},\eqno(**)$$
where for any $A\in\Hom(W,W)$ we denote by $A^\#$
its adjoint operator with respect to $\omega$.
\endproclaim

\demo{Proof}
Let $W_0=\Ker\,\omega|_{\Im A}$,  let
$W_1$ be the orthogonal complement to $W_0$ in $\Im A$ w.r.t.\ the
Hermitian form.
Let $U_2=\Ker A$. Choose a subspace $\tilde U\subset U$
complement to $U_2$.
Then $A$ defines a bijective linear map $\tilde A:\,\tilde U\to\Im A$.
Let $U_0=\tilde A^{-1}(W_0)$, $U_1=\tilde A^{-1}(W_1)$.
We fix a Hermitian form on $U$ such that $U_0$, $U_1$, and $U_2$
are pairwise orthogonal and
claim that the system of equations $(*,**)$ has a solution.

Let $W_2=I\overline{W_0}$, where bar denotes the complex conjugation.
Then $W_0\cap W_2=\{0\}$, the restriction of $\omega$ on $W_0\oplus W_2$
is non-degenerate, and $W_0\oplus W_2$ is orthogonal to $W_1$
both w.r.t.\ $\omega$ and w.r.t.\ the Hermitian form.
Let $W_3$ be the orthogonal complement to $W_0\oplus W_1\oplus W_2$
w.r.t.\ the Hermitian form.

We define the operator $B$ as follows:
$$B|_{W_0}=\tilde A^{-1},\ 
B|_{W_1}=2\tilde A^{-1},\ 
B|_{W_2}=B|_{W_3}=0.$$
Then we have
$$AB|_{W_0}=\Id,\ 
AB|_{W_1}=2\cdot\Id,\ 
AB|_{W_2}=AB|_{W_3}=0,$$
therefore
$$(AB)^\#|_{W_0}=0,\ 
(AB)^\#|_{W_1}=2\cdot\Id,\ 
(AB)^\#|_{W_2}=\Id,\
(AB)^\#|_{W_3}=0.$$
Since $W_0$, $W_1$, $W_2$, and $W_3$ are pairwise 
orthogonal w.r.t.\ the Hermitian
form, it follows that $AB-(AB)^\#$ is a Hermitian operator.
It is easy to check that $(*)$ holds.
Now, we have
$$BA|_{U_0}=\Id,\ BA|_{U_1}=2\cdot\Id,\ BA|_{U_2}=0.$$
Since $U_0$, $U_1$, and $U_2$ are pairwise orthogonal
w.r.t.\ the Hermitian form, it follows that
$BA$ is also a Hermitian operator.
\qed\enddemo

Now we prove Theorem~2.
Until the end of this section ``an algebraic subvariety'' means
``a union of locally closed subvarieties''.
To avoid repetition, we consider the coadjoint case only,
the adjoint case is absolutely similar.

Let $e\in\n$. We need to prove that  
there exists $p\in P$  
such that $\Ad(p)e$ admits a multiple characteristic
$\{h_\chi\}$ with reductive 
$\langle h_\chi\rangle_{alg}$.
At this point, without loss of generality, 
we may assume that $k$ is embedded into $\C$.
By Theorem~1, there exists $p'\in P(\C)$ such that
$\Ad(p')e$ admits a Hermitian multiple characteristic
$\{h'_\chi\}$ (defined over $\C$). It easily follows (see~\cite{Te1})
that $\langle h'_\chi\rangle_{alg}$ is reductive.
Let $W\subset\n$ be a subset of all points that 
admit a multiple characteristic $\{h_\chi\}$ with reductive 
$\langle h_\chi\rangle_{alg}$. We are going to show that
$W$ is a subvariety defined over $k$. Then the argument above
implies that $\Ad(P)e\cap W$ is non-empty over $\C$
hence non-empty over $k$ hence the Theorem.

Let $n$ be the number of positive weights. Consider
the variety of $\ssl_2$-triples 
$$S\subset \mathop\oplus_{\chi<0}\g_{\chi}\oplus\g_0^n\oplus
\mathop\oplus_{\chi>0}\g_{\chi},$$ 
$$S=\bigl\{f_\chi\in\g_{-\chi},
\ h_\chi\in\g_0,\ e_\chi\in\g_\chi\,\big |\,
[e_\chi,f_\chi]=h_\chi,\ [h_\chi,e_\chi]=2e_\chi,
\ [h_\chi,f_\chi]=-2f_\chi\bigr\}.$$
Let $\pi_1$ denote the projection of $S$ on $\g_0^n$,
$\pi_2$ denote the projection of $S$ on $\n$.
These projections are defined over $k$.
Let $R$ denote the set of $n$-tuples of points $(x_1,\ldots,x_n)\in\g_0^n$
such that $\langle x_1,\ldots,x_n\rangle_{alg}$ is reductive.
Then, clearly, $W=\pi_2(\pi_1^{-1}(R))$.
Therefore, it suffices to show that $R$ is a subvariety of $\g_0^n$.
By a theorem of Richardson \cite{Ri},
$\langle x_1,\ldots,x_n\rangle_{alg}$ is reductive if and only if
the orbit of $(x_1,\ldots,x_n)\in\g_0^n$ is closed w.r.t.~the diagonal
action of $L=G_0$.
It remains to notice that if a reductive group $L$ acts on the affine
variety $X$ then the set of closed orbits $Y\subset X$ is a subvariety.
Indeed, one can argue using induction on $\dim X$
and two following observations.
If the action of $G$ on $X$ is not stable (i.e.~a generic $G$-orbit
in $X$ is not closed) then there exists a proper closed subvariety
$X_0\subset X$ that contains all closed orbits, \cite{Vi}.
If the action of $G$ on $X$ is stable, then we can choose
an open $G$-invariant subset $U\in X$  such that all orbits in $U$
are closed in $X$. Then we may pass from $X$ to $X\setminus U$.

\head \S3. Exceptional Groups\endhead

\subhead \S3.0. Comparison Lemma\endsubhead
It follows from Basic Lemma that if all positive weights $\chi\in\ZZ$
are reduced except at most one, then the corresponding parabolic
subgroup is weakly ample.
Therefore, we need to prove Theorem~1 only for parabolic subgroups
such that there exist two or more not reduced $\ZZ$-weights.
Moreover, some not reduced weights may also correspond
to ample components of the grading.
The following lemma provides a lot of examples.

\proclaim{Comparison Lemma}
Let $G'$ and $G''$ be reductive algebraic groups
with Lie algebras $\g'$ and~$\g''$, parabolic subgroups
$P'\subset G'$ and $P''\subset G''$, Levi subgroups $L'\subset P'$ and
$L''\subset P''$. Let $Z'$ and $Z''$ be the centers of $L'$ and $L''$.
Consider positive weights $\chi'\in\ZZ'$ and $\chi''\in\ZZ''$.
Suppose that $[L',L']=[L'',L'']=H$ and $H$-modules
$\g'_{\chi'}$ and $\g''_{\chi''}$ coincide.
Suppose that $H$ is either simple or isomorphic to $\SL_n\times\SL_m$,
in which case $\g'_{\chi'}\simeq\g''_{\chi''}\simeq\C^n\otimes\C^m$
as an $H$-module.
Then, the $\ZZ'$-grading of $\g'$ is ample in degree $\chi'$ iff
the $\ZZ''$-grading of $\g''$ is ample in degree $\chi''$.
\endproclaim

\demo{Proof}
We may assume without loss of generality that $P'$ and $P''$
are maximal parabolic subgroups.
Let $\h$ be the Lie algebra of $H$.
Consider two local Lie algebras
$\g'_{-\chi'}\oplus\g'_0\oplus\g'_{\chi'}$ and
$\g''_{-\chi''}\oplus\g''_0\oplus\g''_{\chi''}$.
We choose $c'\in\z'$ and $c''\in\z''$ such that
$\ad(c')|_{\g'_{\pm\chi'}}=\pm\Id$ and
$\ad(c'')|_{\g''_{\pm\chi''}}=\pm\Id$. 
Then $\g'_0=\C c'\oplus\h$ and $\g''_0=\C c''\oplus\h$.
We fix compact real forms $\k'\subset\g'_0$
and $\k''\subset\g''_0$. Then $\k'=i\R c'\oplus\m$ and
$\k''=i\R c''\oplus\m$, where $\m$ is a compact real form in $\h$.
Suppose that $e'\in \g'_{\chi'}$ and the $\ZZ'$-grading of $\g'$ 
is ample in degree $\chi'$.
Then there exists $f'\in \g'_{-\chi'}$
such that $\langle e',h'=[e',f'],f'\rangle$ is an $\ssl_2$-triple
and $h'\in i\k'$. Let $h'=h'_1+h'_2$, where $h'_1\in i\R c'$ and
$h'_2\in i\m$. 
Let $\h=\oplus\h_p$ be the direct sum of simple Lie algebras.
Then we have the corresponding decomposition $\m=\oplus\m_p$ and
$h'_2=\sum(h'_2)_p$.
We can identify $\g'_0$ and $\g''_0$ as Lie algebras
by assigning $c''$ to $c'$. Then we can identify $\g'_{\chi'}$ with
$\g''_{\chi''}$ and $\g'_{-\chi'}$ with $\g''_{\chi''}$ as their modules.
Let $e''\in\g''_{\chi''}$ corresponds to $e'$ and $f''\in\g''_{-\chi''}$
corresponds to $f'$. Let $h''=[e'',f'']=h''_1+h''_2$,
where $h''_1\in\C c''$ and $h''_2\in\h$, $h''_2=\sum(h''_2)_p$.
We claim that $h''\in i\k''$ and its real multiple is a characterictic of $e''$.
It would follow that
the $\ZZ''$-grading of $\g''$ is ample in degree $\chi''$.

Let $\tilde h'=\tilde h'_1+\tilde h'_2$ be the corresponding
element of $\g'_0$, where $\tilde h'_1\in\C c'$ and $\tilde h'_2\in\h$,
$\tilde h'_2=\sum(\tilde h'_2)_p$.
It is sufficient to prove that $\tilde h'_1=\alpha h'_1$ and
$(\tilde h'_2)_p=\beta_p (h'_2)_p$, 
where $\alpha$ and $\beta$ are positive 
real numbers (here we use our restrictions on $H$).
The commutator maps
$$\Phi':\,\g'_{-\chi'}\otimes\g'_{\chi'}\to \g'_0\quad\text{and}\quad
\Phi'':\,\g''_{-\chi''}\otimes\g''_{\chi''}\to \g''_0$$
can be decomposed as
$\Phi'=\Omega'\circ\Psi'$, 
$\Phi''=\Omega''\circ\Psi''$, where $\Psi'$ and $\Psi''$
are moment maps and
$\Omega':\,(\g_0')^*\to \g_0'$,
$\Omega'':\,(\g_0'')^*\to \g_0''$ are pairings
given $K'$ and $K''$, where $K'$, $K''$ are
restrictions of the Killing forms in $\g'$ and $\g''$
on $\g'_0$ and $\g''_0$.
Under the identifications above we have $\Psi'=\Psi''$.
It remains to notice that $\C c'$ is orthogonal to $\h$
under $K'$, $\C c''$ is orthogonal to $\h$ under $K''$,
restrictions of $K'$ and $K''$ on $(\h)_p$ coincide with the Killing
form of $(\h)_p$ multiplied by a positive real number,
$K'(c',c')$ and $K''(c'',c'')$ are both real and positive.
\qed\enddemo

We can use this lemma and 
examples of ample maximal parabolic subgroups in classical groups
obtained in \cite{Te1} to verify that all parabolic subgroups
not listed in the table at the end of the paper are weakly ample
(this is simple combinatorics).

In particular, we see that all parabolic subgroups in $G_2$, 
$F_4$ and $E_6$ are weakly ample.
Theorem~1 for the entries 1--37 from the table will be checked
in the next sections.
Entries 38--59 will be verified in the sequel of this paper.

\subhead 3.1. Non-degenerate Deformations\endsubhead
A projective variety $X\subset\P(V)$ is called a $2$-variety
if it is cut out by quadrics, i.e.~its homogeneous ideal $I\in\C[V]$
is generated by $I_2$.
For example, if $G$ is a semi-simple Lie group and 
$V$ is an irreducible $G$-module then the projectivization of the 
cone of highest weight vectors in $V$ is a $2$-variety by~\cite{Li}.
Suppose that $X\subset\P(V)$ is a $2$-variety.
A linear map $A:\,\C^2\to V$ is called degenerate
if $\dim\Im A=2$ and $\P(\Im A)\cap X$ is a point.

\proclaim{Proposition 1}\newline
{\rm (A)}
If $A:\,\C^2\to V$ is degenerate,
$B\in\Hom(\C^2,\C^p)$ is non trivial, then there exists
$C\in\Hom(\C^p,V)$ such that $A+CB$ is non-degenerate.\newline
{\rm (B)}
If $A,B:\,\C^2\to V$ are degenerate, then there exists
$E\in\Hom(\C^2,\C^2)$ such that $A+BE$ is non-degenerate.\newline
{\rm (C)}
Suppose $A:\,\C^2\to V$ is degenerate, 
$v\in V$ does not belong to the cone over $X$.
Then there exists $f\in(\C^2)^*$ such that
$A+v\cdot f$ is non-degenerate.
\endproclaim

\demo{Proof}
(A) Clearly, we can choose $C$ such that $\dim\Im(A+CB)<2$.

(B) If $(\Im A)\cap(\Im B)\ne 0$ then there exists
$E\in\Hom(\C^2,\C^2)$ such that we have inequality $\dim \Im(A+BE)<2$ and 
we are done. Otherwise, let $V_0=(\Im A)+(\Im B)$, $X_0=X\cap \P(V_0)$.
Clearly, $\dim V_0=4$ and $X_0$ is a $2$-variety in $\P(V_0)$.
We take skew lines $l_0=\P(\Im A)$ and $l_1=\P(\Im B)$.
Each of them intersects $X_0$ at a point.
For any $E\in\Hom(\C^2,\C^2)$, $\dim\Im(A+BE)=2$.
A line $l\subset\P(V_0)$ is equal to $\P(\Im(A+BE))$ for some
$E\in\Hom(\C^2,\C^2)$ if and only if $l\cap l_1=\emptyset$.

Clearly, $\dim X_0<3$. 
If $\dim X_0<2$ then there exists a line
$l\subset\P(V_0)$ such that $l\cap l_1=l\cap X_0=\emptyset$ and
we are done.
Otherwise, since $l_1\not\subset X_0$,
there exist hyperplanes $H_1,H_2\subset\P(V)$ 
and points $p_1\in H_1$, $p_2\in H_2$ such that
$l_1=H_1\cap H_2$, $p_1,p_2\not\in l_1$, $p_1,p_2\in X_0$.
Then the line $l$ connecting $p_1$ and $p_2$ does not intersect $l_1$
and either intersects $X_0$ in two points or belongs to $X_0$.

(C) Indeed, if $v\in\Im A$, then we can find $f$ such that
$\dim\Im(A+v\cdot f)<2$. Suppose that $v\not\in\Im A$.
Then any $2$-dimensional subspace in $U=\langle v\rangle+\Im A$
not containing $v$ can be realised as $\Im (A+v\cdot f)$ for some $f$.
Let $Y=\P U\cap X$,
$b\in \P U$ be the point corresponding to $v$. 
Then $b\not\in Y$.
We need to find a line $l$ in $\P U$ such that $b\not\in l$
and $l\cap Y$ is not a point.
If $\dim Y=0$, then we can find a line that does not intersect $b$ and $Y$ at all.
If $\dim Y=1$ and $Y$ has a line as an irreducible component,
then we can take this component as a line we are looking for.
Finally, in remaining cases a generic line in $\P U$ intersects $Y$
in more then one (actually, two) points.
\qed\enddemo

\subhead 3.2. Kac Theorem\endsubhead
Consider the $\ZZ$-graded reductive Lie algebra 
$\g=\oplus_{\chi\in\ZZ}\g_{\chi}$, $\chi\ne0$.
Let $\O$ denote the open $L$-orbit in $\g_{\chi}$, $e\in \O$.
Kac has proved in \cite{Ka} that there exists an $\ssl_2$-triple
$\langle e,h,f\rangle$ with $h\in\z$ if and only if the complement
of $\O$ has codimension~1.
In particular, we have the following proposition

\proclaim{Proposition 2}
If $\O$ is open and its complement has codimension~1, then $\O$ is ample.
\endproclaim

\noindent
The list of representations of $[L,L]$ on $\g_{\chi}$ with the above property
can be found in \cite{loc. cit.}.

\subhead 3.3. Witt Theorem\endsubhead
The following proposition is well known as Witt Theorem.

\proclaim{Proposition 3}
Consider a vector space $W$ endowed with 
a non-degenerate quadratic form~$Q$.
Orbits of $\GL(V)\times \OOO(W)$ on $\Hom(V,W)$ are parametrized
by pairs of integers $(i,j)$ such that
$$0\le i\le \min(\dim V,\dim W),\quad
0\le j\le \min(i, \dim W-i).$$
Here $A\in\Hom(V,W)$ corresponds
to the pair $(\dim A, \dim\Ker\, Q|_{\dim A})$.
\endproclaim

With respect to the action of $\GL(V)\times\SO(W)$,
some of these orbits may split into two orbits.
This action appears as the action of a Levi subgroup
of a maximal parabolic subgroup in $G=\SO_n$ on $\g_1$.
It was shown in \cite{Te1} that ample orbits correspond to pairs
$(k,0)$ and $(k,k)$.
More precisely, a linear map $A:\,\C^k\to\C^n$
(or the corresponding tensor $\tilde A\in(\C^k)^*\otimes\C^n$), where
$\C^n$ is endowed with a non-degenerate scalar product,
is called ample if the restriction
of the scalar product in $\C^n$ on $\Im A$ is either non-degenerate
or trivial.

\proclaim{Proposition 4}
Let $V=\C^n$ be endowed with a non-degenerate scalar product.\newline
{\rm (A)} Let $k\le 3$ or $k=n=4$.
Suppose that a linear map $A:\,\C^k\to V$
is not ample, and a linear map
$B:\,\C^k\to \C^p$ is not trivial. 
Then $A+CB$ is ample for some
$C\in\Hom(\C^p, V)$.\newline
{\rm (B)} Let $n\le 4$.
Suppose that the linear map $A:\,\C^k\to V$ is not ample,
$v\in V$ is not trivial.
Then there exists a linear function $f:\,(\C^k)^*$
such that $A+v\cdot f$ is ample.
\endproclaim

\demo{Proof} Simple calculation.\qed\enddemo

The tensor $\tilde A\in\C^k\otimes\C^2\otimes\C^2$ may be considered
as a linear map $A_0:\,(\C^k)^*\to(\C^2\otimes\C^2)$.
The space $\C^2\otimes\C^2$ has a canonical quadratic form $\det$.
$\tilde A$ is called not ample if restriction of $\det$ on $\Im A_0$
is degenerate but not trivial. $\tilde A$ may also be viewed as a linear map
$A:\,(\C^2)^*\to(\C^2\otimes\C^k)$.
Then $\tilde A$ is not ample if and only if $\Im A$ is $2$-dimensional
and $(\Im A)\cap R$ is a line, where $R\subset\C^2\otimes \C^k$
is the variety of rank $1$ matrices.

\proclaim{Proposition 5}\newline
{\rm (A)} Suppose that linear maps $A:\,\C^2\to(\C^2\otimes\C^3)$ and 
$B:\,\C^2\to(\C^2\otimes\Lambda^2\C^3)$ are not ample.
Then there exists $E\in\End(\C^2)\otimes\C^3$ such that
$B+E\circ A$ is ample, here we embed $\C^3$ in $\Hom(\C^3,\Lambda^2\C^3)$
in the obvious way.\newline
{\rm (B)} Suppose that the linear map $A:\,\Lambda^2\C^3\to(\C^2\otimes\C^2)$ is not ample,
$B\in\Hom(\C^3,\C^2)$ is not trivial.
Then there exists a linear map $C\in\Hom(\C^3,\C^2)$
such that $A+B\wedge C$ is ample.\newline
{\rm (C)} Let $v\in\C^3$, $v\ne0$,
$B\in(\Lambda^2\C^3)\otimes\C^2$. 
Then there exists $A\in\C^3\otimes\C^2$ such that
$B+v\wedge A\in(\Lambda^2\C^3)\otimes\C^2$ has rank $2$.
\endproclaim

\demo{Proof}
(A) Direct calculation shows that 
there exists $E_0\in\End(\C^2)\otimes\C^3$ such that
$E_0\circ A$ belongs to the open orbit
in $\Hom(\C^2,\C^2\otimes\Lambda^2\C^3)$ w.r.t.~the group
$G=\GL_2\times\GL_2\times\GL_3$. This complement of this orbit
has codimension~$1$, therefore there exists $\lambda\in\C$ such that
$B+\lambda E_0\circ A$ belongs to this open orbit.
Now we can take $E=\lambda E_0$ by Proposition~2.

(B) Suppose first that $\dim\Im B=2$. We take a basis $\{e_1,e_2,e_3\}$ of $\C^3$
such that $B(e_3)=0$. If we take $C$ such that $C(e_3)=0$,
then $B\wedge C(e_1\wedge e_3)=B\wedge C(e_2\wedge e_3)=0$ 
and $B\wedge C(e_1\wedge e_2)$ can be made arbitrary. Thus we are done
by Proposition~4.A.

Suppose now that $\dim\Im B=1$. Then the linear space of all
possible maps of the form $B\wedge C$ consists of all linear maps
$X:\,\Lambda^2\C^3\to(\C^2\otimes\C^2)$ such that $\Ker X\supset R$
and $\Im X\subset L$, where $R\subset \C^3$ is a fixed
$1$-dimensional subspace and $L\subset \C^2\otimes\C^2$ is a fixed
$2$-dimensional subspace isotropic w.r.t.~$\det$.
We need to consider two cases: $\dim\Im A=2$ and $\dim\Im A=3$.

Let $\dim\Im A=2$. We fix a basis $\{e_1,e_2,e_3\}$ of $\Lambda^2\C^3$
such that $A(e_1)=0$ and $v=A(e_2)$ spans the kernel of restriction
of $\det$ on $\langle v,u\rangle$, where $u=A(e_3)$.
Clearly, we may assume that $R$ is spanned by $e_1$, $e_2$, or $e_3$.
If $R=\C e_1$ and $L\cap\langle v,u\rangle=0$ then we can choose $X$
such that $\Im(A+X)$ is an arbitrary $2$-dimensional subspace not
intersecting $L$ and the claim obviously follows.
If $R=\C e_1$ and $L\cap\langle v,u\rangle\ne0$ then we can choose $X$
such that $\dim(A+X)<2$ and we are done. 
If $R=\C e_2$ and $v\not\in L$ then we can find $w\in L$ such that
$v$ and $w$ are not orthogonal. Then we take the map $X$ such that
$X(e_1)=w$, $X(e_3)=0$. $\Im (A+X)$ is $3$-dimensional and non-degenerate
w.r.t.~$\det$.
If $R=\C e_2$ and $v\in L$ then we take $X$ such that $X(e_1)=0$ and
$u+X(e_3)$ is isotropic w.r.t.~$\det$. This is possible, because
the orthogonal complement to $L$ is $L$ itself, but $u\not\in L$. 
Then $\Im(A+X)$ is $2$-dimensional and isotropic w.r.t.~$\det$.
If $R=\C e_3$ and $v\not\in L$ then we can find $w\in L$ such that
$v$ and $w$ are not orthogonal. Then we take the map $X$ such that
$X(e_1)=w$, $X(e_2)=0$. $\Im (A+X)$ is $3$-dimensional and non-degenerate
w.r.t.~$\det$.
If $R=\C e_3$ and $v\in L$ then we take $X$ such that $X(e_2)=-v$ and
$X(e_1)$ is not orthogonal to $u$ w.r.t.~$\det$. 
This is possible, because
the orthogonal complement to $L$ is $L$ itself, but $u\not\in L$. 
Then $\Im(A+X)$ is $2$-dimensional and non-degenerate w.r.t.~$\det$.

Now let $\dim\Im A=3$.
We fix a basis $\{e_1,e_2,e_3\}$ of $\Lambda^2\C^3$
such that $v=A(e_1)$ spans the kernel of restriction of $\det$
on $\Im A$, $u=A(e_2)$ and $w=A(e_3)$ are isotropic.
We may assume that $R$ either belongs to $\langle e_2,e_3\rangle$
or coincides with $\C e_1$.

Let $R\subset\langle e_2,e_3\rangle$.
If $v\not\in L$, then we take $x\in L$ such that $v$
is not orthogonal to $x$. We take an operator $X$ such that
$X(e_2)=X(e_3)=0$, $X(e_1)=\lambda x$, $\lambda\in\C$.
Then $\Im(A+X)$ is $3$-dimensional and non-degenerate w.r.t.~$\det$
for generic $\lambda$.
If $v\in L$, then we take $X$ such that $X(e_1)=-v$, $X(e_2)=X(e_3)=0$.
Then $\Im (A+X)$ is $2$-dimensional and non-degenerate w.r.t.~$\det$.

Let $R=\C e_1$. If $v\not\in L$, then we take $x\in L$ such that
$v$ is not orthogonal to $x$. We take an operator $X$ such that
$X(e_1)=0$, $X(e_2)=X(e_3)=x$. Then $\Im(A+X)$ is $3$-dimensional
and non-degenerate w.r.t.~$\det$.
If $v\in L$ then either $L=\langle v,u\rangle$ or
$L=\langle v,w\rangle$. It is sufficient to consider only the first case.
We take $X$ such that $X(e_1)=0$, $X(e_2)=-u$, $X(e_3)=0$.
Then $\Im (A+X)=\langle v,w\rangle$ is $2$-dimensional and isotropic.

(C) We fix a basis $\{e_1,e_2,e_3\}$ of $\C^3$ such that $v=e_1$
and a basis $\{f_1,f_2\}$ of $\C^2$.
We take $A=e_2\otimes f_1+e_3\otimes f_2$.
Then $v\wedge A=(e_1\wedge e_2)\otimes f_1+(e_1\wedge e_3)\otimes f_2$
has rank $2$.
Therefore, $B+v\wedge \lambda A$ has rank $2$ for some $\lambda\in\C$.
\qed\enddemo

A tensor 
$\tilde A\in\C^k\otimes\Lambda^2\C^4$
is called ample, if restriction of the quadratic form~$\Pf$
on $\Im A$ is trivial or non-degenerate, where
$A:\,(\C^k)^*\to\Lambda^2\C^4$
is the corresponding map.

\proclaim{Proposition 6}\newline
{\rm (A)} Suppose that linear maps 
$A:\,\C^2\to\Lambda^2\C^4$ and $B:\,\C^2\otimes\C^2\to\Lambda^3\C^4$
are not ample. Then there exists a map $C:\,\C^2\to\C^4$ such that
$B+A\wedge C$ is ample.\newline
{\rm (B)} Suppose that linear maps $A:\,\Lambda^2\C^4\to\C^2$ 
and $B:\,\C^4\to(\C^2\otimes\C^2)$
are not ample. Then there exists a map $C:\,\C^4\to(\C^2)^*\otimes\C$ such that
$A+B\wedge C$ is ample.\newline
{\rm (C)} Suppose that a linear map $A:\,\C^2\to\Lambda^2\C^4$ 
is not ample, $w\in\C^4$ is not trivial. 
Then there exists a map $B:\,\C^2\to\C^4$ such that
$A+w\wedge B$ is ample.\newline
{\rm(D)} Suppose that $A\in\C^4\otimes\Lambda^2\C^4$ is not ample.
If $B\in\Lambda^2\C^4\otimes\Lambda^3\C^4$ is not ample, then
there exists $C\in\C^4\otimes\C^4$ such that
$B+A\wedge C$ is ample.\newline
{\rm(E)} Suppose that $A\in\C^4\otimes\Lambda^2\C^4$ is not ample.
If $u\in\C^4$, $u\ne0$, then there exists
$C=\sum v_i\otimes w_i\in\C^4\otimes\C^4$ such that
$A+\sum v_i\otimes (w_i\wedge u)$ is~ample. 
\endproclaim

\demo{Proof}
(A) It suffices to find $C$ such that $A\wedge C$ belongs to the open orbit,
because it is ample.
We fix bases $\{e_1,e_2,e_3,e_4\}$ of $\C^4$, $\{f_1,f_2\}$ of the first $\C^2$,
$\{g_1,g_2\}$ of the second $\C^2$. We may assume
that $A(f_1)=e_1\wedge e_2+e_3\wedge e_4$, $A(f_2)=e_1\wedge e_3$.
Consider $C$ such that $C(g_1)=e_4$, $C(g_2)=e_2$.
Then $A\wedge C$ is an isomorphism.

(B) If $\dim B=3$, then there exists a map $C:\,\C^4\to(\C^2)^*\otimes\C$ such that
$A+B\wedge C$ belongs to the open orbit, hence it is ample.
Let $\dim B=2$.
Simple calculations show that we need to prove the following claim.

Let $\{e_1,e_2,e_3,e_4\}$ (resp.~$\{f_1,f_2\}$)
be a basis of $\C^4$ (resp.~of $\C^2$).
Suppose  that the map $A:\,\C^2\to\Lambda^2\C^4$ is not ample.
Then we claim that there exist vectors $v,u\in\C^4$ such that
the map $B$ is ample, where $B(f_1)=A(f_1)+e_1\wedge u$,
$B(f_2)=A(f_2)+e_2\wedge u+e_1\wedge v$.
Indeed, $\dim \Im A=2$ and we may assume that either $A(f_1)$ or $A(f_2)$
spans the kernel of restriction of $\Pf$ on $\Im A$.
Consider the first case. If $A(f_1)=e_1\wedge w$, then we take $v=0$, $u=-w$.
If $A(f_1)$ is not of the form $e_1\wedge w$, then there exists $v$ such that
$A(f_1)\wedge e_1\wedge v\ne0$. If we take $u=0$, then $B$ is ample.
Consider the second case. Since $A(f_1)\wedge A(f_1)\ne0$,
we can find $w$ such that $A(f_1)\wedge e_1\wedge w\ne0$.
If $v=\lambda w$ and $u=0$, then $B$ is ample for generic $\lambda\in\C$.

(C)
We fix a basis $\{e_1,e_2,e_3,e_4\}$ of $\C^4$ such that
$w=e_1$ and a basis $\{f_1,f_2\}$ such that
$A(f_1)$ spans the kernel of restriction of $\Pf$ on
$\Im A$. We need to prove that there exist $u,v\in\C^4$
such that the linear map $C:\,\C^2\to\Lambda^2\C^4$ is ample,
where $C(f_1)=A(f_1)+e_1\wedge u$, $C(f_2)=A(f_2)+e_1\wedge v$.
If $A(f_1)=e_1\wedge x$, then it suffices to take $u=-x$.
Otherwise, we take $u=0$ and we choose $v$ such that
$A(f_1)\wedge e_1\wedge v\ne0$. Then $\dim\Im C=2$ and restriction
of $\Pf$ on $\Im C$ is non-degenerate.

(D) and (E).
Let $\tilde A:\,(\C^4)^*\to\Lambda^2\C^4$
be the corresponding map. Since restriction of $\Pf$
on $\Im\tilde A$ is degenerate but not-trivial,
in the suitable bases $A$ has one of the following forms:
$$A_1=e_1\otimes (f_1\wedge f_2+f_3\wedge f_4)+e_2\otimes(f_1\wedge f_3),$$
$$A_2=e_1\otimes (f_1\wedge f_2+f_3\wedge f_4)+
e_2\otimes(f_1\wedge f_3)+e_3\otimes(f_1\wedge f_4),$$
$$A_3=e_1\otimes (f_1\wedge f_2+f_3\wedge f_4)+
e_2\otimes(f_1\wedge f_3)+e_3\otimes(f_3\wedge f_4),$$
$$A_4=e_1\otimes (f_1\wedge f_2+f_3\wedge f_4)+
e_2\otimes(f_1\wedge f_3)+e_3\otimes(f_1\wedge f_4)+e_4\otimes(f_1\wedge f_2),$$
$$A_5=e_1\otimes (f_1\wedge f_2+f_3\wedge f_4)+
e_2\otimes(f_1\wedge f_3)+e_3\otimes(f_1\wedge f_4)+e_4\otimes(f_2\wedge f_4).$$

(D) We take
$C=e_3\otimes f_2+e_4\otimes f_4$.
Then
$$A\wedge C=(e_1\wedge e_3)\otimes (f_2\wedge f_3\wedge f_4)-
(e_2\wedge e_3)\otimes (f_1\wedge f_2\wedge f_3)+$$
$$+(e_1\wedge e_4)\otimes (f_1\wedge f_2\wedge f_4)+
(e_2\wedge e_4)\otimes (f_1\wedge f_3\wedge f_4).$$
Restriction of $\Pf$ on the $4$-dimensional image of the corresponding map
is non-degenerate, so $A\wedge C$ belongs to the open orbit, which is ample
by the results of \S3.2. Therefore, $B+A\wedge(\lambda C)$ is also ample
for some $\lambda\in \C$.

(E) Let $A=A_1$. If $u\in\langle f_1,f_3\rangle$, then 
we take $C=e_2\otimes u'$ such that $u'\wedge u=-f_1\wedge f_3$.
If $u\not\in\langle f_1,f_3\rangle$, then there exists $u'$ such that
$f_1\wedge f_3\wedge u'\wedge u\ne0$ and we take $C=e_1\otimes u'$.

Let $A=A_2$ or $A_4$. 
If $u\not\in\langle f_1,f_3\rangle$, then 
there exists $u'$ such that
$f_1\wedge f_3\wedge u'\wedge u\ne0$ and we take $C=e_3\otimes u'$.
If $u\not\in\langle f_1,f_4\rangle$, then 
there exists $u'$ such that
$f_1\wedge f_4\wedge u'\wedge u\ne0$ and we take $C=e_2\otimes u'$.
Finally, if $u=\lambda f_1$, then we take
$C=(e_2\otimes f_3+e_3\otimes f_4)/\lambda$.

Let $A=A_3$. If $u\in\langle f_1,f_3\rangle$, then 
we take $C=e_2\otimes u'$ such that $u'\wedge u=-f_1\wedge f_3$.
If $u\not\in\langle f_1,f_3\rangle$, then there exists $u'$ such that
$f_1\wedge f_3\wedge u'\wedge u\ne0$ and we take $C=e_3\otimes u'$.

Let $A=A_5$. If $u\in\langle f_1,f_4\rangle$, then 
we take $C=e_3\otimes u'$ such that $u'\wedge u=-f_1\wedge f_4$.
If $u\not\in\langle f_1,f_4\rangle$, then there exists $u'$ such that
$f_1\wedge f_4\wedge u'\wedge u\ne0$ and we take $C=e_1\otimes u'$.
\qed\enddemo

\subhead 3.4. Spinors\endsubhead

We recall the definition of half-spinor representations.
Let $V=\C^{2m}$ be an even-dimensional vector space with
a non-degenerate symmetric scalar product $(\cdot,\cdot)$.
Let $\Cl(V)$ be the Clifford algebra of $V$.
Recall that $\Cl(V)$ is in fact a superalgebra,
$\Cl(V)=\Cl^0(V)\oplus\Cl^1(V)$, where
$\Cl^0(V)$ (resp.~$\Cl^1(V)$) is a linear span of elements
of the form $v_1\cdot\ldots\cdot v_r$, $v_i\in V$, $r$ is even 
(resp.~$r$ is odd). Then we have
$$\Spin(V)=\{a\in\Cl^0(V)\,|\,a=v_1\cdot\ldots\cdot v_r,\ v_i\in V,\ Q(v_i,v_i)=1\}.$$
$\Spin(V)$ is a connected simply-connected
algebraic group.

Given $a\in\Cl(V)$, $a=v_1\cdot\ldots\cdot v_r$, $v_i\in V$,
let $\overline a=(-1)^rv_r\cdot\ldots\cdot v_1$.
This is a well-defined involution of $\Cl(V)$.
Using it, we may define an action $R$ 
of $\Spin(V)$ on $V$ by formula $R(a)v=a\cdot v\cdot\overline a$.
Then $R$ is a double covering $\Spin(V)\to\SO(V)$.

Let $U\subset V$ be a maximal isotropic subspace, hence $\dim U=m$.
Take also any maximal isotropic subspace $U'$
such that $U\oplus U'=V$. For any $v\in U$ we define
an operator $\rho(v)\in\End(\Lambda^*U)$ by the formula
$$\rho(v)\cdot u_1\wedge\ldots\wedge u_r=v\wedge u_1\wedge\ldots u_r.$$
For any $v\in U'$ we define an 
operator $\rho(v)\in\End(\Lambda^*U)$ by the formula
$$\rho(v)\cdot u_1\wedge\ldots\wedge u_r=\sum_{i=1}^r(-1)^{i-1}(v,u_i)
\wedge u_1\wedge\ldots u_{i-1}\wedge u_{i+1}\wedge\ldots\wedge u_r.$$
These operators are well-defined and
therefore by linearity we have a linear map
$\rho:\, V\to\End(\Lambda^*U)$.
It is easy to check that for any $v\in V$ we have
$\rho(v)^2=(v,v)\Id$.
Therefore we have a homomorphism
of associative algebras
$\Cl(V)\to\End(\Lambda^*U)$, i.e., $\Lambda^*U$
is a $\Cl(V)$-module, easily seen to be irreducible.
Therefore, $\Lambda^*U$ is also a $\Spin(V)$-module 
(spinor representation), now reducible.
However, the even part $\Lambda^{ev}U$ is an irreducible
$\Spin(V)$-module, called the half-spinor module $S^+$.
Another irreducible half-spinor module $S^-$ is defined as $\Lambda^{od}U$.

The map $\rho$ defines a morphism of $\Spin(V)$ modules
$V\otimes S^{\pm}\to S^{\mp}$. 
If $m$ is odd, then $S^+$ is dual to $S^-$.
If $m$ is even, 
then the spinor representation $\Lambda^*U$ admits an interesting
non-degenerate bilinear form $(\cdot,\cdot)$ defined as follows.
Let $\det\in\Lambda^mU$ be a fixed non-trivial element.
Then $(u,v)$ is equal to the coefficient at $\det$ of 
the element $(-1)^{\left[{\deg u\over 2}\right]}u\wedge v$. 
Obviously, $S^+$ is orthogonal to $S^-$ and restriction of 
$(u,v)$ on $S^\pm$ is orthogonal if $m=4k$ and symplectic
if $m=4k+2$. In particular, $S^+$ and $S^-$ are self-dual if $m$ is even.

If $m=4$, then $V$, $S^+$ and $S^-$ are twisted forms of each other
w.r.t.~outer isomorphisms of $\Spin_{8}$ (triality principle).
Let $V=\C^8$ be a vector space equipped with a non-degenerate
scalar product. We denote by $S^+$, $S^-$ the corresponding
half-spinor modules. If $R$ denotes $V$, $S^+$, or $S^-$, then
a tensor $A\in\C^k\otimes R$ is called ample
if restriction of the scalar product of $R$ on $\Im\tilde A$
is trivial or non-degenerate, where $\tilde A:(\C^k)^*\to R$
is the corresponding map.

\proclaim{Proposition 7}\newline
{\rm (A)}
Let $k=2,3$.
If $A=\sum a_i\otimes s_i\in \C^k\otimes S^+$ and 
$B\in\C^k\otimes S^-$ are not ample,
then there exists $v\in V$ such that 
$C=B+\sum a_i\otimes \rho(v)s_i$ is ample.\newline
{\rm (B)}
If $A=\sum a_i\otimes s_i\in \C^3\otimes S^+$ and 
$B\in\Lambda^2\C^3\otimes S^-$ are not ample,
then there exists $x=\sum v_i\otimes s_i\in \C^3\otimes V$ 
such that $C=B+\sum_{i,j} (a_i\wedge v_j)\otimes \rho(x_j)s_i$ is ample.\newline
{\rm (C)}
If $A\in \C^3\otimes S^-$ is not ample and $s\in S^+$ is not trivial,
then there exists $B=\sum a_i\otimes w_i\in \C^3\otimes V$ 
such that $C=A+\sum a_i\otimes \rho(w_i)s$ is ample.
\endproclaim

\demo{Proof}
(A) Let $k=2$.
We choose a basis $f_1,f_2$ of $\C^2$ and $e_1,e_2,e_3,e_4$ of $U$
(maximal isotropic subspace of $V$) such that
$$A=f_1\otimes 1+f_2\otimes (e_1\wedge e_2+e_3\wedge e_4),\quad
B=f_1\otimes x_1+f_2\otimes x_2.$$
If $x_1\in U\subset\Lambda^{od}U$ then we take $v=-x_1$.
Otherwise, there exists $v'\in U'$ such that 
$\rho(v')(e_1\wedge e_2+e_3\wedge e_4)$ is not orthogonal to $x_1$.
Then we take $v=\lambda v'$: for generic $\lambda$ restriction
of $(\cdot,\cdot)$ on $\Im\tilde C$ is non-degenerate,
hence $C$ is ample.

Let $k=3$. 
We choose a basis $\{f_1,f_2,f_3\}$ of $\C^3$ 
and $\{e_1,e_2,e_3,e_4\}$ of $U$
(maximal isotropic subspace of~$V$) such that $A$ has one
of the following forms:
$$A_1=f_1\otimes (e_1\wedge e_2+e_3\wedge e_4)+f_2\otimes 1,$$
$$A_2=f_1\otimes 1+f_2\otimes (e_1\wedge e_2)+f_3\otimes (e_3\wedge e_4),$$
$$A_3=f_1\otimes (1+e_1\wedge e_2\wedge e_3\wedge e_4)+
f_2\otimes e_1\wedge e_2+f_3\otimes e_1\wedge e_3.$$
Let
$$B=f_1\otimes x_1+f_2\otimes x_2+f_3\otimes x_3.$$

{\bf Let $\bold{A=A_1}$.} 
We need to show, that there exist $u,v\in U$ such that the following element if ample:
$$C=f_1\otimes\bigl(x_1+u+v\wedge (e_1\wedge e_2+e_3\wedge e_4)\bigr)
+f_2\otimes (x_2+v)+f_3\otimes x_3.$$
If $x_3\not\in U$, then there exists $v'\in U$ such that
$( v,x_3)\ne0$ and $u'\in U$ such that 
$(e_1\wedge e_2+e_3\wedge e_4)\wedge u'\wedge v'\ne0$.
We take $v=\lambda v'$, $u=\lambda u'$.
Then for generic $\lambda$ the image of $\tilde C$ is $3$-dimensional
and non-degenerate, hence $C$ is ample.

Suppose now that $x_3\in U$.
If $x_3=0$, then we finish the proof as in the case $k=2$.

Let $x_3\ne0$.
If $x_2\not\in U$, then there exist $v'\in U$ such that
$( x_3, v'\wedge (e_1\wedge e_2+e_3\wedge e_4))\ne0$,
$x_2\wedge v'\ne0$. We take $u=0$, $v=\lambda v'$.
Then for generic $\lambda$ the image of $\tilde C$ is $3$-dimensional
and non-degenerate, hence $C$ is ample.

Suppose now that $x_2\in U$. Then we choose $u$ and $v$ such that
$u+v\wedge (e_1\wedge e_2+e_3\wedge e_4)=-x_1$.
Then the image of $\tilde C$ is isotropic.

{\bf Let $\bold{A=A_2}$.} 
We need to show, that there exist $v\in U$, $u\in\langle e_1,e_2\rangle$,
$w\in\langle e_3,e_4\rangle$ such that the following element if ample:
$$C=f_1\otimes(x_1+v)+f_2\otimes (x_2+v\wedge(e_1\wedge e_2)+u)
+f_3\otimes (x_3+v\wedge(e_3\wedge e_4)+w).$$
 
If $x_1\not\in U$, then there exists $v'\in U$ such that
$( v,x_1)\ne0$ and 
$u\in\langle e_1,e_2\rangle$ (or
$w\in\langle e_3,e_4\rangle$)
such that $u\wedge v\wedge e_3\wedge e_4\ne0$
(or $w\wedge v\wedge e_1\wedge e_2\ne0$)
We take $v=\lambda v'$, $u=\lambda u'$, $w=0$ (or 
$w=\lambda w'$, $u=0$).
Then for generic $\lambda$ the image of $\tilde C$ is $3$-dimensional
and non-degenerate, hence $C$ is ample.

Suppose now that $x_1\in U$.
Then we may suppose that $x_1=0$ after taking $v=-x_1$.
Now we are going to find $u$ and $w$.

If $x_3$ is not perpendicular to $\langle e_1,e_2\rangle$
(resp.~$x_2$ is not perpendicular to $\langle e_3,e_4\rangle$),
then we take $u'$ such that $(u',x_3)\ne0$, $w'=0$ 
(resp.~take $w'$ such that $(w',x_2)\ne0$, $u'=0$)
and set $u=\lambda u'$, $w=\lambda w'$. Then $\Im\tilde C$
is $2$-dimensional and non-degenerate, hence ample, for generic $\lambda$.

Let $x_3\perp\langle e_1,e_2\rangle$, $x_2\perp\langle e_3,e_4\rangle$.

If $(x_3,x_3)=0$, then $(x_2,x_3)=0$, otherwise $B$ is ample.
Then $(x_2,x_2)\ne0$, otherwise $B$ is ample.
Hence $x_2$ is not perpendicular to $\langle e_1,e_2\rangle$
and we take $w=0$ and $v$ such that $(x_2,x_2)+2(x_2,u)=0$

Suppose, therefore, that $(x_3,x_3)\ne0$, and, similarly, that
$(x_2,x_2)\ne0$. It follows that
$x_2$ is not perpendicular to $\langle e_1,e_2\rangle$, 
$x_3$ is not perpendicular to $\langle e_3,e_4\rangle$.
We take $u'$, $w'$ such that
$(x_2,u)\ne0$, $(x_3,w')\ne0$ We set $u=\lambda u'$, $w=\lambda w'$.
Then $\Im\tilde C$ is $2$-dimensional and non-degenerate for generic $\lambda$,
hence $C$ is ample for these values of $\lambda$.

{\bf Let $\bold{A=A_3}$.}
There exists $x\in V$ such that 
$\rho(x)(1+e_1\wedge e_2\wedge e_3\wedge e_4)=-x_1$.
Therefore, without loss of generality we may assume that $x_1=0$.
If $B$ is ample, then there is nothing to prove. Otherwise,
the restriction of a scalar form $R$ on $\langle x_2,x_3\rangle$
has one-dimensional kernel. After a suitable change of bases we may 
assume that this kernel is spanned by $x_2$.
There exists $x'\in V$ such that 
$\rho(x')(1+e_1\wedge e_2\wedge e_3\wedge e_4)$ is not orthogonal to $x_2$.
We take $x=\lambda x'$. Then for generic $\lambda$
the image of $\tilde C$ is three-dimensional and non-degenerate
w.r.t.~$R$.

(B)
We choose a basis $\{f_1,f_2,f_3\}$ of $\C^2$ 
and $\{e_1,e_2,e_3,e_4\}$ of $U$
(maximal isotropic subspace of~$V$) such that $A$ has one
of the following forms:
$$A_1=f_1\otimes (e_1\wedge e_2+e_3\wedge e_4)+f_2\otimes 1,$$
$$A_2=f_1\otimes 1+f_2\otimes (e_1\wedge e_2)+f_3\otimes (e_3\wedge e_4),$$
$$A_3=f_1\otimes (1+e_1\wedge e_2\wedge e_3\wedge e_4)+
f_2\otimes e_1\wedge e_2+f_3\otimes e_1\wedge e_3.$$
Let 
$$B=(f_2\wedge f_3)\otimes x_1+(f_1\wedge f_3)\otimes x_2+(f_1\wedge f_2)\otimes x_3.$$
We take an isotropic subspace $U'\subset V$ complement to $U$ and choose
a basis
$\{e_1^*,e_2^*,e_3^*,e_4^*\}$ of $U'$ such that $(e_i,e_j^*)=\delta_{ij}$.

{\bf Let $\bold{A=A_3}$.} 
We take 
$$x=\lambda(f_3\otimes e_1+f_1\otimes e_4+f_2\otimes e_1^*+f_3\otimes e_3^*).$$
Then the image of $\tilde D$ is $3$-dimensional and non-degenerate,
where
$D=\sum_{i,j} (a_i\wedge v_j)\otimes \rho(x_j)s_i$.
Therefore, for generic $\lambda$, the image of $\tilde C$ 
is $3$-dimensional and non-degenerate, hence $C$ is ample.

{\bf Let $\bold{A=A_2}$.} 
The same proof as above, but for
$$x=\lambda(f_3\otimes e_1^*+f_2\otimes e_3^*+f_1\otimes 
(e_1+e_2+e_3+e_1^*)).$$

{\bf Let $\bold{A=A_1}$.} 
The same proof as above, but for
$$x=\lambda(f_3\otimes e_1+f_3\otimes e_1^*+f_2\otimes e_2).$$

(C)
If $(s,s)\ne0$, then $\rho(V)s=S^-$,
otherwise,  $\rho(V)s=U_0$, where $U_0\subset S^-$ is
a maximal isotropic subspace. It is clear now that
we can find $B$ such that $\Im \tilde C$ is isotropic.
\qed\enddemo

\head \S4. $E_7$ and $E_8$\,: The Zoo\endhead

We associate to each parabolic subgroup its coloured
Dynkin diagram. To each graded component $\g_{\chi}$, $\chi>0$
we associate the coloured Dynkin diagram of the corresponding
parabolic subgroup with vertices indexed as follows.
The lattice $\ZZ$ is isomorphic to a sublattice in the root lattice
spanned by simple roots marked white on the coloured Dynkin diagram.
We write the corresponding coefficient of $\chi$ near each white vertex.
The representation of $[\l,\l]$ on $\g_{\chi}$ is irreducible and
the Dynkin diagram of $[\l,\l]$ is equal to the black subdiagram 
of the coloured Dynkin diagram. We write the numerical labels of $\chi$
near each black vertice (omitting zeroes).
We shall pick several positive weights  
and call them {\it twisting weights}.
All reduced positive weights of the form $\chi+\sum\pm\mu_i$, where $\chi>0$
is not reduced and all $\mu_i>0$ are twisting, are called
{\it rubbish weights}.
We shall draw a diagram with the set of vertices given
by not reduced and rubbish weights and arrows indexed 
by twisting weights: an arrow $\mu$ has a tail $\chi_1$ and a head $\chi_2$
if and only if $[\g_{\chi_1},\g_{\mu}]=\g_{\chi_2}$.
For any $e\in\n$ {\rm(}resp.~$e\in\g/\p^-${\rm)}
we shall try to find an element $u\in U$ {\rm(}resp.~$u\in U^-${\rm)}
such that all not reduced graded components of the element $\Ad(u)e$ 
belong to ample orbits, except at most one.
Here $U\subset P$ is the subgroup with a Lie algebra generated
by all twisting weights components,
$U^-\subset P^-$ is the subgroup with a Lie algebra generated
by all graded components $\g_{\chi}$, where $-\chi$ is twisting.
The standard strategy will be to decrease the number of non-reduced
components by applying elements of the form $\exp p$,
where $p$ belongs to the graded component of some twisting weight.
These `elementary transformations' will be made in the
special order, because they may change some other components as well.
We shall use big Latin letters for not reduced weights,
small Latin letters for rubbish weights, and numbers for twisting weights.
To save space, we shall use abbreviations of the form
$$L\LRA{n}{X.Y}M.$$ 
This means
that if $x_L\subset\g_L$ is not ample (or not trivial if $L$ is
a reduced weight), then we can apply the element $\exp(p_n)$,
where $p_n\in\g_n$, to make $x_M\in\g_M$ ample, and
reason for this is given in Proposition~X.Y.
Parabolic subgroups from the table will be joined in several groups,
in each group the proof is similar.

\subhead Case 1A\endsubhead
This case includes parabolic subgroups  1 and 2. 
We shall consider the 1st. 

\medskip
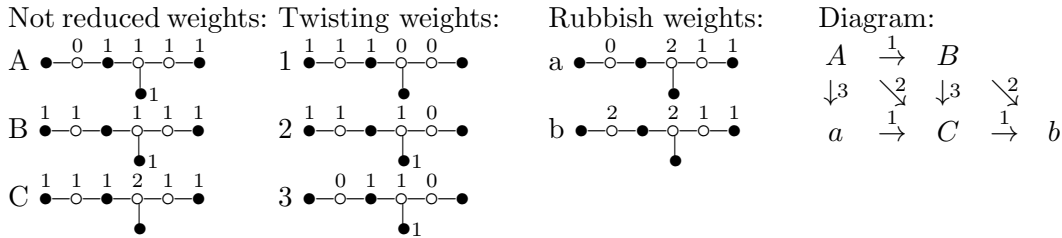
\hbox{{\vtop{\hsize=0.25\hsize\parindent0pt
Not reduced weights:

A $\Ess{\B^{}}{\W^{0}}{\B^{1}}{\W^{1}}{\W^{1}}{\B^{1}}\B{1}$

B $\Ess{\B^{1}}{\W^{1}}{\B^{}}{\W^{1}}{\W^{1}}{\B^{1}}\B{1}$

C $\Ess{\B^{1}}{\W^{1}}{\B^{1}}{\W^{2}}{\W^{1}}{\B^{1}}\B{}$

}}{
\vtop{\hsize=0.25\hsize\parindent0pt
Twisting weights:

1 $\Ess{\B^{1}}{\W^{1}}{\B^{1}}{\W^{0}}{\W^{0}}{\B^{}}\B{}$

2 $\Ess{\B^{1}}{\W^{1}}{\B^{}}{\W^{1}}{\W^{0}}{\B^{}}\B{1}$

3 $\Ess{\B^{}}{\W^{0}}{\B^{1}}{\W^{1}}{\W^{0}}{\B^{}}\B{1}$

}}{
\vtop{\hsize=0.25\hsize\parindent0pt
Rubbish weights:

a $\Ess{\B^{}}{\W^{0}}{\B^{}}{\W^{2}}{\W^{1}}{\B^{1}}\B{}$

b $\Ess{\B^{}}{\W^{2}}{\B^{}}{\W^{2}}{\W^{1}}{\B^{1}}\B{}$

}}{
\vtop{\hsize=0.25\hsize\parindent0pt
Diagram:

$\matrix
A&\ra1&B&&\cr
\da3&\sea2&\da3&\sea2&\cr
a&\ra1&C&\ra1&b\cr
\endmatrix$
}}}
\medskip

Consider the coadjoint case, $x\in\n$.
Then 
$A\LRA{1}{1.B}B$, $A\LRA{2}{1.B}C$. If $x_A$ is ample, then
$B\LRA{3}{1.B}C$.

Adjoint case, $x\in\g/\p^-$.
$b\LRA{-2}{1.A}B$, $b\LRA{-1}{1.A}C$.
Let $x_b=0$. Then $C\LRA{-3}{1.B}B$, $C\LRA{-2}{1.B}A$.
If $x_C$ is ample, then $B\LRA{-1}{1.B}A$.

\subhead Case 1B\endsubhead
This case includes parabolic subgroups 3, 4, 5, 6, 7, 8, 9.
The proof is similar to the previous case (using Proposition 1.A and 1.B)
but combinatorics is more involved. To save space, we omit this proof
and refer the reader to the even more complicated Case~5.

\subhead Case 2A\endsubhead 
Parabolic subgroup number 19.

\hbox{{\vtop{\hsize=0.25\hsize\parindent0pt

A $\Ess{\B^{}}{\W^{0}}{\B^{}}{\B^{1}}{\W^{1}}{\B^{1}}\B{}$

B $\Ess{\B^{1}}{\W^{1}}{\B^{1}}{\B^{}}{\W^{1}}{\B^{1}}\B{}$

}}{
\vtop{\hsize=0.25\hsize\parindent0pt

1 $\Ess{\B^{1}}{\W^{1}}{\B^{}}{\B^{}}{\W^{0}}{\B^{}}\B{1}$

}}{
\vtop{\hsize=0.25\hsize\parindent0pt

a $\Ess{\B^{}}{\W^{2}}{\B^{}}{\B^{}}{\W^{1}}{\B^{1}}\B{}$

}}{
\vtop{\hsize=0.25\hsize\parindent0pt
${}$ 

$A\ra{1}B\ra{1}a$

}}}
\medskip

\noindent
Adjoint and coadjoint cases follow from Propositions~6.A, 6.B and~1.A.

\subhead Case 2B\endsubhead 
This case includes parabolic subgroups with numbers 20 and 21.
We shall consider 20.

\hbox{{\vtop{\hsize=0.25\hsize\parindent0pt
A $\Ess{\B^{1}}{\B^{}}{\B^{}}{\W^{1}}{\W^{1}}{\B^{1}}\B{1}$

}}{
\vtop{\hsize=0.25\hsize\parindent0pt
B $\Ess{\B^{}}{\B^{1}}{\B^{}}{\W^{2}}{\W^{1}}{\B^{1}}\B{}$

}}{
\vtop{\hsize=0.25\hsize\parindent0pt
1 $\Ess{\B^{1}}{\B^{}}{\B^{}}{\W^{1}}{\W^{0}}{\B^{}}\B{1}$

}}{
\vtop{\hsize=0.25\hsize\parindent0pt
${}$

$A\ra{1}B$

}}}
\medskip

\noindent
Both adjoint and coadjoint cases follow from Proposition~6.A and
Proposition~6.B.

\subhead Case 2C\endsubhead 
This case includes parabolic subgroups 22, 23, 24, 25, 26, 27. 
We shall consider the 26.

\hbox{{\vtop{\hsize=0.25\hsize\parindent0pt
A $\Ees{\B^{1}}{\B^{}}{\B^{}}{\W^{1}}{\B^{1}}{\W^{1}}{\B^{1}}\W{1}$

B $\Ees{\B^{}}{\B^{1}}{\B^{}}{\W^{2}}{\B^{}}{\W^{1}}{\B^{1}}\W{1}$

C $\Ees{\B^{}}{\B^{1}}{\B^{}}{\W^{2}}{\B^{1}}{\W^{2}}{\B^{}}\W{1}$

}}{
\vtop{\hsize=0.25\hsize\parindent0pt
1 $\Ees{\B^{1}}{\B^{}}{\B^{}}{\W^{1}}{\B^{1}}{\W^{0}}{\B^{}}\W{0}$

2 $\Ees{\B^{1}}{\B^{}}{\B^{}}{\W^{1}}{\B^{}}{\W^{1}}{\B^{1}}\W{0}$

3 $\Ees{\B^{}}{\B^{}}{\B^{}}{\W^{0}}{\B^{1}}{\W^{1}}{\B^{1}}\W{0}$

}}{
\vtop{\hsize=0.25\hsize\parindent0pt
a $\Ees{\B^{1}}{\B^{}}{\B^{}}{\W^{1}}{\B^{}}{\W^{2}}{\B^{}}\W{1}$

b $\Ees{\B^{}}{\B^{}}{\B^{1}}{\W^{3}}{\B^{}}{\W^{2}}{\B^{}}\W{1}$

}}{
\vtop{\hsize=0.25\hsize\parindent0pt
${}$ 

$\matrix
A&\ra1&B&&\cr
\da3&\sea2&\da3&\sea2&\cr
a&\ra1&C&\ra1&b\cr
\endmatrix$

}}}
\medskip

Coadjoint case, $x\in\n$.
Then $A\LRA{1}{6.B}B$, $A\LRA{2}{6.B}C$.
If $x_A$ is ample, then
$B\LRA{3}{1.B}C$.

Adjoint case, $x\in\g/\p^-$.
Then $b\LRA{-2}{6.C}B$ and $b\LRA{-1}{6.C}C$.
Let $x_b=0$. Then $C\LRA{-3}{1.B}B$ and
$C\LRA{-2}{6.A}A$.
Let $x_C$ be ample.
Then $B\LRA{-1}{6.A}A$.

\subhead Case 2D\endsubhead 
This case includes parabolic subgroups 27, 28, 29. 
We shall consider 29.

\hbox{{\vtop{\hsize=0.25\hsize\parindent0pt
A $\Ees{\B^{}}{\W^{0}}{\B^{}}{\B^{1}}{\B^{}}{\W^{1}}{\B^{1}}\W{1}$

B $\Ees{\B^{1}}{\W^{1}}{\B^{1}}{\B^{}}{\B^{}}{\W^{1}}{\B^{1}}\W{1}$

C $\Ees{\B^{1}}{\W^{1}}{\B^{}}{\B^{1}}{\B^{}}{\W^{2}}{\B^{}}\W{1}$

}}{
\vtop{\hsize=0.25\hsize\parindent0pt
1 $\Ees{\B^{1}}{\W^{1}}{\B^{}}{\B^{}}{\B^{1}}{\W^{0}}{\B^{}}\W{0}$

2 $\Ees{\B^{1}}{\W^{1}}{\B^{}}{\B^{}}{\B^{}}{\W^{1}}{\B^{1}}\W{0}$

3 $\Ees{\B^{}}{\W^{0}}{\B^{1}}{\B^{}}{\B^{}}{\W^{1}}{\B^{1}}\W{0}$

}}{
\vtop{\hsize=0.25\hsize\parindent0pt
a $\Ees{\B^{}}{\W^{0}}{\B^{}}{\B^{}}{\B^{1}}{\W^{2}}{\B^{}}\W{1}$

b $\Ees{\B^{}}{\W^{2}}{\B^{}}{\B^{}}{\B^{}}{\W^{1}}{\B^{1}}\W{1}$

c $\Ees{\B^{}}{\W^{2}}{\B^{1}}{\B^{}}{\B^{}}{\W^{2}}{\B^{}}\W{1}$

}}{
\vtop{\hsize=0.25\hsize\parindent0pt
${}$ 

$\matrix
A&\ra1&B&\ra1&b\cr
\da3&\sea2&\da3&\sea2&\da3\cr
a&\ra1&C&\ra1&c\cr
\endmatrix$

}}}
\medskip

Coadjoint case, $x\in\n$.
Then $A\LRA{1}{6.A}B$,
$A\LRA{2}{1.B}C$.
If $x_A$ is ample, then $B\LRA{3}{6.B}C$.

Adjoint case, $x\in\g/\p^-$.
Then $c\LRA{-1}{6.C}C$ and $c\LRA{-2}{4.A}B$.
Let $x_c=0$.
Then $C\LRA{-3}{6.A}B$ 
and $C\LRA{-2}{1.B}A$.
Let $x_C$ be ample.
If $x_b=0$, then $B\LRA{-1}{6.B}A$.
If $x_b\ne0$ then $b\LRA{-1}{1.A}B$.

\subhead Case 2E\endsubhead 
This case includes parabolic subgroups 30, 31. 
We shall consider 31.

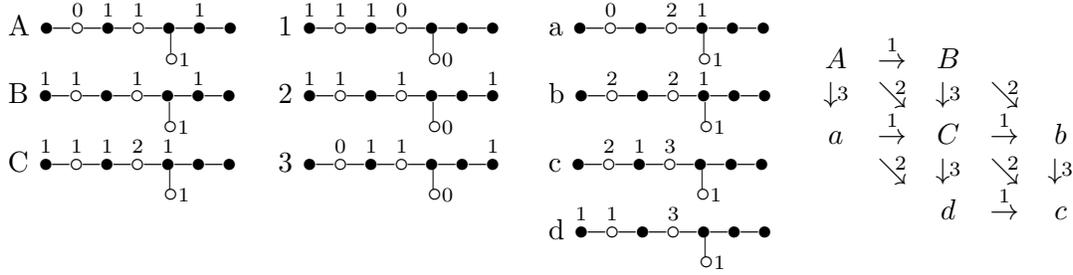
\hbox{{\vtop{\hsize=0.25\hsize\parindent0pt
A $\Ees{\B^{}}{\W^{0}}{\B^{1}}{\W^{1}}{\B^{}}{\B^{1}}{\B^{}}\W{1}$

B $\Ees{\B^{1}}{\W^{1}}{\B^{}}{\W^{1}}{\B^{}}{\B^{1}}{\B^{}}\W{1}$

C $\Ees{\B^{1}}{\W^{1}}{\B^{1}}{\W^{2}}{\B^{1}}{\B^{}}{\B^{}}\W{1}$

}}{
\vtop{\hsize=0.25\hsize\parindent0pt
1 $\Ees{\B^{1}}{\W^{1}}{\B^{1}}{\W^{0}}{\B^{}}{\B^{}}{\B^{}}\W{0}$

2 $\Ees{\B^{1}}{\W^{1}}{\B^{}}{\W^{1}}{\B^{}}{\B^{}}{\B^{1}}\W{0}$

3 $\Ees{\B^{}}{\W^{0}}{\B^{1}}{\W^{1}}{\B^{}}{\B^{}}{\B^{1}}\W{0}$

}}{
\vtop{\hsize=0.25\hsize\parindent0pt
a $\Ees{\B^{}}{\W^{0}}{\B^{}}{\W^{2}}{\B^{1}}{\B^{}}{\B^{}}\W{1}$

b $\Ees{\B^{}}{\W^{2}}{\B^{}}{\W^{2}}{\B^{1}}{\B^{}}{\B^{}}\W{1}$

c $\Ees{\B^{}}{\W^{2}}{\B^{1}}{\W^{3}}{\B^{}}{\B^{}}{\B^{}}\W{1}$

d $\Ees{\B^{1}}{\W^{1}}{\B^{}}{\W^{3}}{\B^{}}{\B^{}}{\B^{}}\W{1}$

}}{
\vtop{\hsize=0.25\hsize\parindent0pt
${}$ 

$\matrix
A&\ra1&B&&\cr
\da3&\sea2&\da3&\sea2&\cr
a&\ra1&C&\ra1&b\cr
&\sea2&\da3&\sea2&\da3\cr
&&d&\ra1&c\cr
\endmatrix$

}}}

Coadjoint case, $x\in\n$.
Then $A\LRA{1}{1.B}B$ and
$A\LRA{2}{6.A}C$.
If $x_A$ is ample, then $B\LRA{3}{6.A}C$.

Adjoint case, $x\in\g/\p^-$.
If $x_c\ne0$ or $x_b\ne0$, then we first 
apply $\exp(p_{-3})$ to make $x_b\ne0$ if necessary.
Then $b\LRA{-2}{6.C}B$ and
$b\LRA{-1}{4.A}C$.
If $x_c,x_b=0$ then  
$d\LRA{-3}{1.A}C$ and $B\LRA{-1}{1.B}A$.
If $x_c,x_b,x_d=0$ then $C\LRA{-3}{6.B}B$ and
$C\LRA{-2}{6.B}A$. If $x_C$ is ample, then
$B\LRA{-1}{1.B}A$.

\subhead Case 3\endsubhead 
This case includes parabolic subgroups 32, 33. 
We shall consider 33.

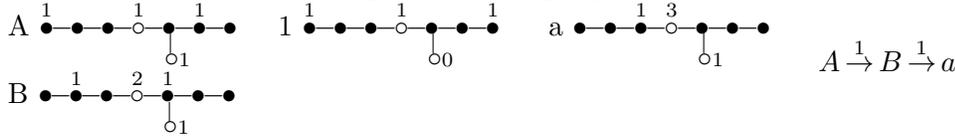
\hbox{{\vtop{\hsize=0.25\hsize\parindent0pt
A $\Ees{\B^{1}}{\B^{}}{\B^{}}{\W^{1}}{\B^{}}{\B^{1}}{\B^{}}\W{1}$

B $\Ees{\B^{}}{\B^{1}}{\B^{}}{\W^{2}}{\B^{1}}{\B^{}}{\B^{}}\W{1}$

}}{
\vtop{\hsize=0.25\hsize\parindent0pt
1 $\Ees{\B^{1}}{\B^{}}{\B^{}}{\W^{1}}{\B^{}}{\B^{}}{\B^{1}}\W{0}$

}}{
\vtop{\hsize=0.25\hsize\parindent0pt
a $\Ees{\B^{}}{\B^{}}{\B^{1}}{\W^{3}}{\B^{}}{\B^{}}{\B^{}}\W{1}$

}}{
\vtop{\hsize=0.25\hsize\parindent0pt
${}$ 

$A\ra{1}B\ra{1}a$

}}}
\medskip

\noindent
Adjoint and coadjoint cases
follow from Proposition~6.D and Proposition~6.E.

\subhead Case 4A\endsubhead 
In this case we study the parabolic subgroup number 34.

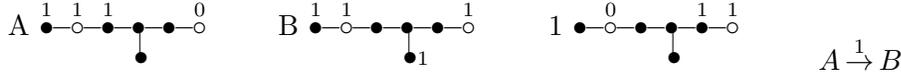
\hbox{{\vtop{\hsize=0.25\hsize\parindent0pt
A $\Ess{\B^{1}}{\W^{1}}{\B^{1}}{\B^{}}{\B^{}}{\W^{0}}\B{}$

}}{
\vtop{\hsize=0.25\hsize\parindent0pt
B $\Ess{\B^{1}}{\W^{1}}{\B^{}}{\B^{}}{\B^{}}{\W^{1}}\B{1}$

}}{
\vtop{\hsize=0.25\hsize\parindent0pt
1 $\Ess{\B^{}}{\W^{0}}{\B^{}}{\B^{}}{\B^{1}}{\W^{1}}\B{}$

}}{
\vtop{\hsize=0.25\hsize\parindent0pt
${}$ 

$A\ra{1}B$

}}}
\medskip

\noindent
Both adjoint and coadjoint cases follow from Proposition~7.A. 

\subhead Case 4B\endsubhead
This case includes parabolic subgroups 35, 36.
We shall consider 36.

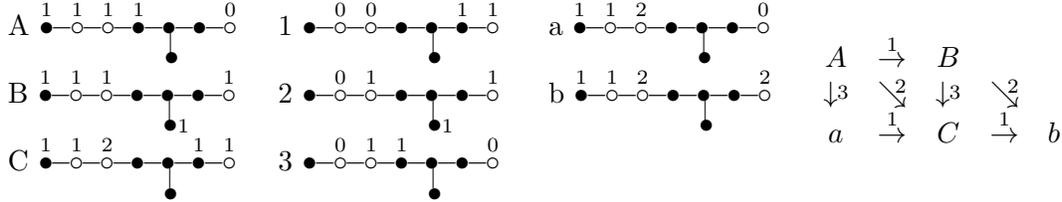
\hbox{{\vtop{\hsize=0.25\hsize\parindent0pt
A $\Ees{\B^{1}}{\W^{1}}{\W^{1}}{\B^{1}}{\B^{}}{\B^{}}{\W^{0}}\B{}$

B $\Ees{\B^{1}}{\W^{1}}{\W^{1}}{\B^{}}{\B^{}}{\B^{}}{\W^{1}}\B{1}$

C $\Ees{\B^{1}}{\W^{1}}{\W^{2}}{\B^{}}{\B^{}}{\B^{1}}{\W^{1}}\B{}$

}}{
\vtop{\hsize=0.25\hsize\parindent0pt
1 $\Ees{\B^{}}{\W^{0}}{\W^{0}}{\B^{}}{\B^{}}{\B^{1}}{\W^{1}}\B{}$

2 $\Ees{\B^{}}{\W^{0}}{\W^{1}}{\B^{}}{\B^{}}{\B^{}}{\W^{1}}\B{1}$

3 $\Ees{\B^{}}{\W^{0}}{\W^{1}}{\B^{1}}{\B^{}}{\B^{}}{\W^{0}}\B{}$

}}{
\vtop{\hsize=0.25\hsize\parindent0pt
a $\Ees{\B^{1}}{\W^{1}}{\W^{2}}{\B^{}}{\B^{}}{\B^{}}{\W^{0}}\B{}$

b $\Ees{\B^{1}}{\W^{1}}{\W^{2}}{\B^{}}{\B^{}}{\B^{}}{\W^{2}}\B{}$

}}{
\vtop{\hsize=0.25\hsize\parindent0pt
${}$

$\matrix
A&\ra1&B&&\cr
\da3&\sea2&\da3&\sea2&\cr
a&\ra1&C&\ra1&b\cr
\endmatrix$

}}}
\medskip

\noindent
Easy digram-walking shows that everything follows from Propositions~7.A and~1.A.

\subhead Case 4C\endsubhead
This case includes the parabolic subgroup 37.

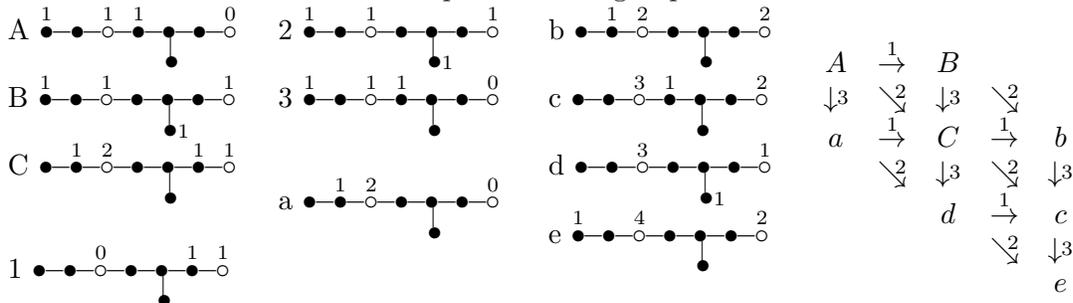
\hbox{{\vtop{\hsize=0.25\hsize\parindent0pt
A $\Ees{\B^{1}}{\B^{}}{\W^{1}}{\B^{1}}{\B^{}}{\B^{}}{\W^{0}}\B{}$

B $\Ees{\B^{1}}{\B^{}}{\W^{1}}{\B^{}}{\B^{}}{\B^{}}{\W^{1}}\B{1}$

C $\Ees{\B^{}}{\B^{1}}{\W^{2}}{\B^{}}{\B^{}}{\B^{1}}{\W^{1}}\B{}$

\bigskip

1 $\Ees{\B^{}}{\B^{}}{\W^{0}}{\B^{}}{\B^{}}{\B^{1}}{\W^{1}}\B{}$

}}{
\vtop{\hsize=0.25\hsize\parindent0pt
2 $\Ees{\B^{1}}{\B^{}}{\W^{1}}{\B^{}}{\B^{}}{\B^{}}{\W^{1}}\B{1}$

3 $\Ees{\B^{1}}{\B^{}}{\W^{1}}{\B^{1}}{\B^{}}{\B^{}}{\W^{0}}\B{}$

\bigskip

a $\Ees{\B^{}}{\B^{1}}{\W^{2}}{\B^{}}{\B^{}}{\B^{}}{\W^{0}}\B{}$

}}{
\vtop{\hsize=0.25\hsize\parindent0pt
b $\Ees{\B^{}}{\B^{1}}{\W^{2}}{\B^{}}{\B^{}}{\B^{}}{\W^{2}}\B{}$

c $\Ees{\B^{}}{\B^{}}{\W^{3}}{\B^{1}}{\B^{}}{\B^{}}{\W^{2}}\B{}$

d $\Ees{\B^{}}{\B^{}}{\W^{3}}{\B^{}}{\B^{}}{\B^{}}{\W^{1}}\B{1}$

e $\Ees{\B^{1}}{\B^{}}{\W^{4}}{\B^{}}{\B^{}}{\B^{}}{\W^{2}}\B{}$

}}{
\vtop{\hsize=0.25\hsize\parindent0pt
${}$ 

$\matrix
A&\ra1&B&&\cr
\da3&\sea2&\da3&\sea2&\cr
a&\ra1&C&\ra1&b\cr
&\sea2&\da3&\sea2&\da3\cr
&&d&\ra1&c\cr
&&&\sea2&\da3\cr
&&&&e\cr
\endmatrix$

}}}

The coadjoint case follows from 
Proposition~7.A and Proposition~7.B.
The adjoint case easily follows from Proposition~7.A, Proposition~7.B, 
Proposition 4.A, and Proposition~7.C.

\subhead Case 5A\endsubhead This case contains
parabolic subgroup 10.

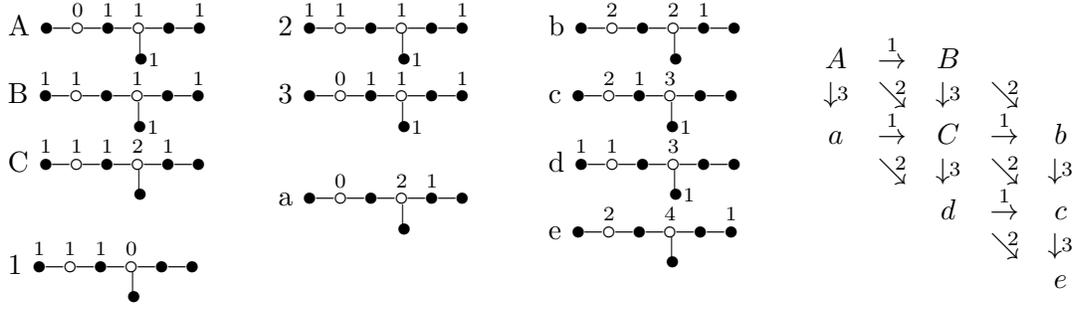
\hbox{{\vtop{\hsize=0.25\hsize\parindent0pt
A $\Ess{\B^{}}{\W^{0}}{\B^{1}}{\W^{1}}{\B^{}}{\B^{1}}\B{1}$

B $\Ess{\B^{1}}{\W^{1}}{\B^{}}{\W^{1}}{\B^{}}{\B^{1}}\B{1}$

C $\Ess{\B^{1}}{\W^{1}}{\B^{1}}{\W^{2}}{\B^{1}}{\B^{}}\B{}$

\bigskip

1 $\Ess{\B^{1}}{\W^{1}}{\B^{1}}{\W^{0}}{\B^{}}{\B^{}}\B{}$

}}{
\vtop{\hsize=0.25\hsize\parindent0pt
2 $\Ess{\B^{1}}{\W^{1}}{\B^{}}{\W^{1}}{\B^{}}{\B^{1}}\B{1}$

3 $\Ess{\B^{}}{\W^{0}}{\B^{1}}{\W^{1}}{\B^{}}{\B^{1}}\B{1}$

\bigskip

a $\Ess{\B^{}}{\W^{0}}{\B^{}}{\W^{2}}{\B^{1}}{\B^{}}\B{}$

}}{
\vtop{\hsize=0.25\hsize\parindent0pt
b $\Ess{\B^{}}{\W^{2}}{\B^{}}{\W^{2}}{\B^{1}}{\B^{}}\B{}$

c $\Ess{\B^{}}{\W^{2}}{\B^{1}}{\W^{3}}{\B^{}}{\B^{}}\B{1}$

d $\Ess{\B^{1}}{\W^{1}}{\B^{}}{\W^{3}}{\B^{}}{\B^{}}\B{1}$

e $\Ess{\B^{}}{\W^{2}}{\B^{}}{\W^{4}}{\B^{}}{\B^{1}}\B{}$

}}{
\vtop{\hsize=0.25\hsize\parindent0pt
${}$

$\matrix
A&\ra1&B&&\cr
\da3&\sea2&\da3&\sea2&\cr
a&\ra1&C&\ra1&b\cr
&\sea2&\da3&\sea2&\da3\cr
&&d&\ra1&c\cr
&&&\sea2&\da3\cr
&&&&e\cr
\endmatrix$

}}}

Consider the coadjoint case, let $x\in\n$.
Then
$A\LRA{1}{1.B}B$,
$A\LRA{2}{5.A}C$.
If $x_A$ is ample, then
$B\LRA{3}{5.A}C$.

Adjoint case, $x\in\g/\p^-$. 
Suppose first that either $x_b$, $x_c$, or $x_e$ is not trivial.
Applying $\exp(p_{-3})$ if necessary, we may assume that $x_b\ne0$.
Then 
$b\LRA{-2}{4.A}B$,
$b\LRA{-1}{4.A}C$.
Let $x_b=x_c=x_e=0$.
Then
$d\LRA{-3}{1.A}C$,
$B\LRA{-1}{1.B}A$.
Let $x_d=0$.
Then 
$C\LRA{-3}{5.A}B$,
$C\LRA{-2}{1.B}A$.
Let $x_C$ be ample, then
$B\LRA{-1}{1.B}A$.

\subhead Case 5B\endsubhead
This case includes parabolic subgroups 11, 12, 13.
We shall consider the 11.

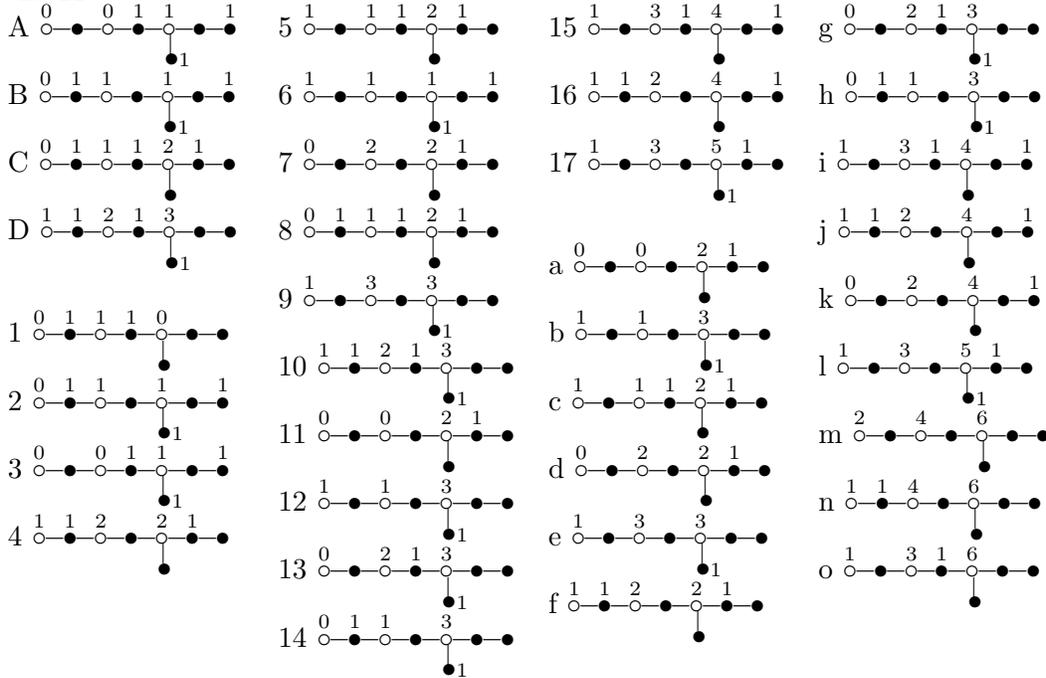
\hbox{{\vtop{\hsize=0.25\hsize\parindent0pt
A $\Ees{\W^{0}}{\B^{}}{\W^{0}}{\B^{1}}{\W^{1}}{\B^{}}{\B^{1}}\B{1}$

B $\Ees{\W^{0}}{\B^{1}}{\W^{1}}{\B^{}}{\W^{1}}{\B^{}}{\B^{1}}\B{1}$

C $\Ees{\W^{0}}{\B^{1}}{\W^{1}}{\B^{1}}{\W^{2}}{\B^{1}}{\B^{}}\B{}$

D $\Ees{\W^{1}}{\B^{1}}{\W^{2}}{\B^{1}}{\W^{3}}{\B^{}}{\B^{}}\B{1}$

\bigskip\noindent

1 $\Ees{\W^{0}}{\B^{1}}{\W^{1}}{\B^{1}}{\W^{0}}{\B^{}}{\B^{}}\B{}$

2 $\Ees{\W^{0}}{\B^{1}}{\W^{1}}{\B^{}}{\W^{1}}{\B^{}}{\B^{1}}\B{1}$

3 $\Ees{\W^{0}}{\B^{}}{\W^{0}}{\B^{1}}{\W^{1}}{\B^{}}{\B^{1}}\B{1}$

4 $\Ees{\W^{1}}{\B^{1}}{\W^{2}}{\B^{}}{\W^{2}}{\B^{1}}{\B^{}}\B{}$

}}{
\vtop{\hsize=0.25\hsize\parindent0pt

5 $\Ees{\W^{1}}{\B^{}}{\W^{1}}{\B^{1}}{\W^{2}}{\B^{1}}{\B^{}}\B{}$

6 $\Ees{\W^{1}}{\B^{}}{\W^{1}}{\B^{}}{\W^{1}}{\B^{}}{\B^{1}}\B{1}$

7 $\Ees{\W^{0}}{\B^{}}{\W^{2}}{\B^{}}{\W^{2}}{\B^{1}}{\B^{}}\B{}$

8 $\Ees{\W^{0}}{\B^{1}}{\W^{1}}{\B^{1}}{\W^{2}}{\B^{1}}{\B^{}}\B{}$

9 $\Ees{\W^{1}}{\B^{}}{\W^{3}}{\B^{}}{\W^{3}}{\B^{}}{\B^{}}\B{1}$

10 $\Ees{\W^{1}}{\B^{1}}{\W^{2}}{\B^{1}}{\W^{3}}{\B^{}}{\B^{}}\B{1}$

11 $\Ees{\W^{0}}{\B^{}}{\W^{0}}{\B^{}}{\W^{2}}{\B^{1}}{\B^{}}\B{}$

12 $\Ees{\W^{1}}{\B^{}}{\W^{1}}{\B^{}}{\W^{3}}{\B^{}}{\B^{}}\B{1}$

13 $\Ees{\W^{0}}{\B^{}}{\W^{2}}{\B^{1}}{\W^{3}}{\B^{}}{\B^{}}\B{1}$

14 $\Ees{\W^{0}}{\B^{1}}{\W^{1}}{\B^{}}{\W^{3}}{\B^{}}{\B^{}}\B{1}$

}}{
\vtop{\hsize=0.25\hsize\parindent0pt

15 $\Ees{\W^{1}}{\B^{}}{\W^{3}}{\B^{1}}{\W^{4}}{\B^{}}{\B^{1}}\B{}$

16 $\Ees{\W^{1}}{\B^{1}}{\W^{2}}{\B^{}}{\W^{4}}{\B^{}}{\B^{1}}\B{}$

17 $\Ees{\W^{1}}{\B^{}}{\W^{3}}{\B^{}}{\W^{5}}{\B^{1}}{\B^{}}\B{1}$

\bigskip\noindent

a $\Ees{\W^{0}}{\B^{}}{\W^{0}}{\B^{}}{\W^{2}}{\B^{1}}{\B^{}}\B{}$

b $\Ees{\W^{1}}{\B^{}}{\W^{1}}{\B^{}}{\W^{3}}{\B^{}}{\B^{}}\B{1}$

c $\Ees{\W^{1}}{\B^{}}{\W^{1}}{\B^{1}}{\W^{2}}{\B^{1}}{\B^{}}\B{}$

d $\Ees{\W^{0}}{\B^{}}{\W^{2}}{\B^{}}{\W^{2}}{\B^{1}}{\B^{}}\B{}$

e $\Ees{\W^{1}}{\B^{}}{\W^{3}}{\B^{}}{\W^{3}}{\B^{}}{\B^{}}\B{1}$

f $\Ees{\W^{1}}{\B^{1}}{\W^{2}}{\B^{}}{\W^{2}}{\B^{1}}{\B^{}}\B{}$

}}{
\vtop{\hsize=0.25\hsize\parindent0pt

g $\Ees{\W^{0}}{\B^{}}{\W^{2}}{\B^{1}}{\W^{3}}{\B^{}}{\B^{}}\B{1}$

h $\Ees{\W^{0}}{\B^{1}}{\W^{1}}{\B^{}}{\W^{3}}{\B^{}}{\B^{}}\B{1}$

i $\Ees{\W^{1}}{\B^{}}{\W^{3}}{\B^{1}}{\W^{4}}{\B^{}}{\B^{1}}\B{}$

j $\Ees{\W^{1}}{\B^{1}}{\W^{2}}{\B^{}}{\W^{4}}{\B^{}}{\B^{1}}\B{}$

k $\Ees{\W^{0}}{\B^{}}{\W^{2}}{\B^{}}{\W^{4}}{\B^{}}{\B^{1}}\B{}$

l $\Ees{\W^{1}}{\B^{}}{\W^{3}}{\B^{}}{\W^{5}}{\B^{1}}{\B^{}}\B{1}$

m $\Ees{\W^{2}}{\B^{}}{\W^{4}}{\B^{}}{\W^{6}}{\B^{}}{\B^{}}\B{}$

n $\Ees{\W^{1}}{\B^{1}}{\W^{4}}{\B^{}}{\W^{6}}{\B^{}}{\B^{}}\B{}$

o $\Ees{\W^{1}}{\B^{}}{\W^{3}}{\B^{1}}{\W^{6}}{\B^{}}{\B^{}}\B{}$

}}}

\noindent
$A\ra{1}B$\quad
$a\ra{1}C\ra{1}d$\quad
$b\ra{1}D\ra{1}e$\quad
$c\ra{1}f$\quad
$h\ra{1}g$\quad
$j\ra{1}i$\quad
$o\ra{1}n$\newline
$A\ra{2}C\ra{2}g$\quad
$B\ra{2}d$\quad
$c\ra{2}D\ra{2}i$\quad
$a\ra{2}h\ra{2}k$\quad
$b\ra{2}j\ra{2}l\ra{2}n$\quad
$f\ra{2}e$\newline
$A\ra{3}a$\quad
$B\ra{3}C\ra{3}h$\quad
$f\ra{3}D\ra{3}j$\quad
$c\ra{3}b$\quad
$d\ra{3}g\ra{3}k$\quad
$e\ra{3}i\ra{3}l\ra{3}o$\newline
$A\ra{4}D$\quad
$B\ra{4}e$\quad
$C\ra{4}i$\quad
$a\ra{4}j\ra{4}m$\quad
$h\ra{4}l$\quad
$k\ra{4}n$\newline
$A\ra{5}b$\quad
$B\ra{5}D$\quad
$C\ra{5}j$\quad
$d\ra{5}i\ra{5}m$\quad
$g\ra{5}l$\quad
$k\ra{5}o$\newline
$A\ra{6}c$\quad
$B\ra{6}f$\quad
$C\ra{6}D$\quad
$a\ra{6}b$\quad
$d\ra{6}e$\quad
$g\ra{6}i$\quad
$h\ra{6}j$\quad
$k\ra{6}l\ra{6}m$\newline
$A\ra{7}g$\quad
$a\ra{7}k$\quad
$b\ra{7}l$\quad
$c\ra{7}i$\quad
$j\ra{7}n$\newline
$A\ra{8}h$\quad
$B\ra{8}g$\quad
$C\ra{8}k$\quad
$D\ra{8}l$\quad
$c\ra{8}j\ra{8}o$\quad
$f\ra{8}i\ra{8}n$\newline
$A\mathop{\rightarrow}\limits^{9}i$\quad
$a\mathop{\rightarrow}\limits^{9}l$\quad
$b\mathop{\rightarrow}\limits^{9}m$\quad
$h\mathop{\rightarrow}\limits^{9}n$\newline
$A\ra{10}j$\quad
$B\ra{10}i$\quad
$C\ra{10}l$\quad
$D\ra{10}m$\quad
$g\ra{10}n$\quad
$h\ra{10}o$\newline
$B\mathop{\rightarrow}\limits^{11}h$\quad
$d\mathop{\rightarrow}\limits^{11}k$\quad
$e\mathop{\rightarrow}\limits^{11}l$\quad
$f\mathop{\rightarrow}\limits^{11}j$\quad
$i\mathop{\rightarrow}\limits^{11}o$\quad\quad
$B\ra{12}j$\quad
$d\ra{12}l$\quad
$e\ra{12}m$\quad
$g\ra{12}o$\newline
$A\ra{13}k$\quad
$D\ra{13}n$\quad
$b\ra{13}o$\quad
$c\ra{13}l$\quad\quad
$B\ra{14}k$\quad
$D\ra{14}o$\quad
$e\ra{14}n$\quad
$f\ra{14}l$\newline
$A\ra{15}l$\quad
$C\ra{15}n$\quad
$a\ra{15}o$\newline
$c\ra{15}m$\quad\quad
$B\ra{16}l$\quad
$C\ra{16}o$\quad
$d\ra{16}n$\quad
$f\ra{16}m$\quad\quad
$A\ra{17}o$\quad
$B\ra{17}n$

\medskip

Consider the coadjoint case, $x\in\n$.
Then $A\LRA{1}{1.B}B$,
$A\LRA{2}{5.A}C$,
$A\LRA{4}{4.B}D$.
Suppose that $x_A$ is ample, then
$B\LRA{3}{5.A}C$,
$B\LRA{5}{4.B}D$.
Let  $x_B$ be ample, then 
$C\LRA{6}{4.B}D$.

Adjoint case, $x\in\g/\p^-$.
Then
$n\LRA{-13}{1.A}D$,
$n\LRA{-15}{1.A}C$,
$n\LRA{-17}{1.A}B$.
Let $x_n=0$, then
$o\LRA{-14}{1.A}D$,
$o\LRA{-16}{1.A}C$,
$o\LRA{-17}{1.A}A$.
Let $x_o=0$, $x_l\ne0$ or $x_m\ne0$ (in the latter case
we apply $\exp(p_{-6}$ if necessary to make $x_l\ne0$).
Then 
$l\LRA{-10}{4.A}C$,
$l\LRA{-15}{5.C}A$,
$l\LRA{-16}{5.C}B$.
Let $x_l=x_m=0$, then
$i\LRA{-2}{1.A}D$,
$i\LRA{-4}{5.C}C$,
$i\LRA{-10}{4.A}B$.
Let $x_i=0$, then
$j\LRA{-3}{1.A}D$,
$j\LRA{-5}{5.C}C$,
$j\LRA{-10}{4.A}A$.
Let $x_j=0$, then
$k\LRA{-8}{4.A}C$,
$k\LRA{-13}{4.A}A$,
$k\LRA{-14}{4.A}B$.
Let $x_k=0$, then
$g\LRA{-2}{1.A}C$,
$g\LRA{-8}{1.A}B$,
$g\LRA{-7}{4.B}A$.
Let $x_g=0$, then
$h\LRA{-3}{1.A}C$,
$h\LRA{-8}{1.A}A$,
$h\LRA{-11}{4.B}B$.
Let $x_h=0$, $x_e\ne0$.
Applying $\exp(p_{-6})$ if necessary, we also make $x_d=0$.
Then 
$e\LRA{-1}{4.A}D$,
$e\LRA{-4}{1.A}B$,
$C\LRA{-2}{5.A}A$,
Let $x_e=0$, then
$D\LRA{-6}{4.B}C$,
$D\LRA{-5}{4.B}B$,
$D\LRA{-4}{4.B}A$.
Let  $x_D$ be ample, then
$d\LRA{-1}{4.A}C$,
$d\LRA{-2}{4.A}B$.
Let $x_d=0$, then
$C\LRA{-2}{5.A}A$,
$C\LRA{-3}{5.A}B$.
Let $x_C$ be ample. Then
$B\LRA{-1}{1.B}A$.

\subhead Case 5C\endsubhead This case includes parabolic subgroups 14, 15.
We shall consider number 14.

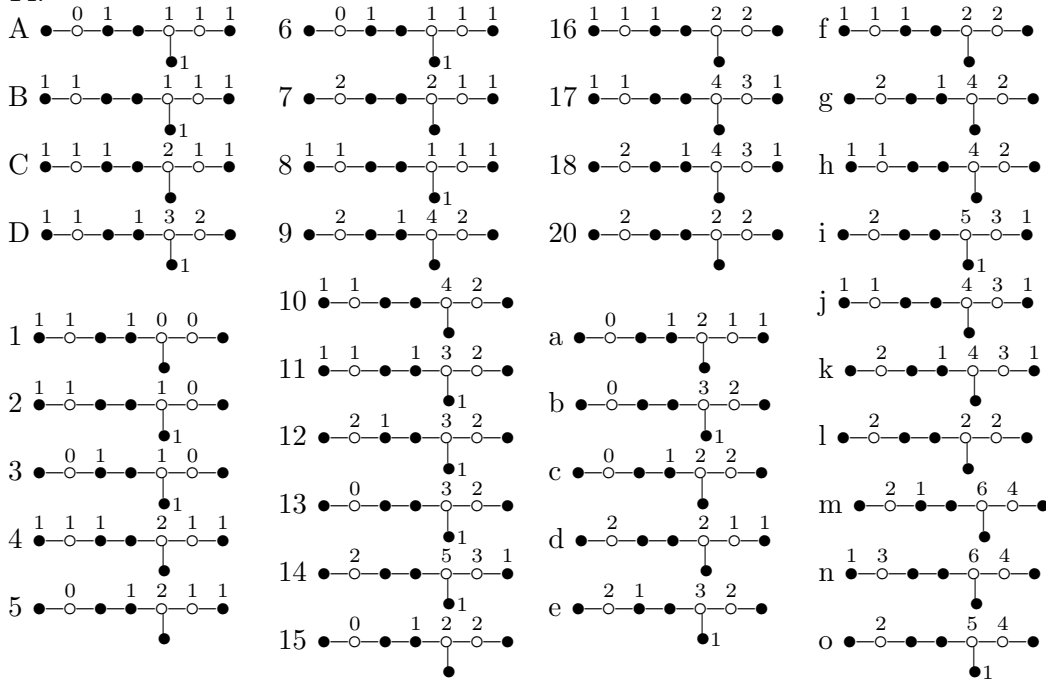
\hbox{{\vtop{\hsize=0.25\hsize\parindent0pt
A $\Ees{\B^{}}{\W^{0}}{\B^{1}}{\B^{}}{\W^{1}}{\W^{1}}{\B^{1}}\B{1}$

B $\Ees{\B^{1}}{\W^{1}}{\B^{}}{\B^{}}{\W^{1}}{\W^{1}}{\B^{1}}\B{1}$

C $\Ees{\B^{1}}{\W^{1}}{\B^{1}}{\B^{}}{\W^{2}}{\W^{1}}{\B^{1}}\B{}$

D $\Ees{\B^{1}}{\W^{1}}{\B^{}}{\B^{1}}{\W^{3}}{\W^{2}}{\B^{}}\B{1}$

\bigskip\noindent

1 $\Ees{\B^{1}}{\W^{1}}{\B^{}}{\B^{1}}{\W^{0}}{\W^{0}}{\B^{}}\B{}$

2 $\Ees{\B^{1}}{\W^{1}}{\B^{}}{\B^{}}{\W^{1}}{\W^{0}}{\B^{}}\B{1}$

3 $\Ees{\B^{}}{\W^{0}}{\B^{1}}{\B^{}}{\W^{1}}{\W^{0}}{\B^{}}\B{1}$

4 $\Ees{\B^{1}}{\W^{1}}{\B^{1}}{\B^{}}{\W^{2}}{\W^{1}}{\B^{1}}\B{}$

5 $\Ees{\B^{}}{\W^{0}}{\B^{}}{\B^{1}}{\W^{2}}{\W^{1}}{\B^{1}}\B{}$

}}{
\vtop{\hsize=0.25\hsize\parindent0pt

6 $\Ees{\B^{}}{\W^{0}}{\B^{1}}{\B^{}}{\W^{1}}{\W^{1}}{\B^{1}}\B{1}$

7 $\Ees{\B^{}}{\W^{2}}{\B^{}}{\B^{}}{\W^{2}}{\W^{1}}{\B^{1}}\B{}$

8 $\Ees{\B^{1}}{\W^{1}}{\B^{}}{\B^{}}{\W^{1}}{\W^{1}}{\B^{1}}\B{1}$

9 $\Ees{\B^{}}{\W^{2}}{\B^{}}{\B^{1}}{\W^{4}}{\W^{2}}{\B^{}}\B{}$

10 $\Ees{\B^{1}}{\W^{1}}{\B^{}}{\B^{}}{\W^{4}}{\W^{2}}{\B^{}}\B{}$

11 $\Ees{\B^{1}}{\W^{1}}{\B^{}}{\B^{1}}{\W^{3}}{\W^{2}}{\B^{}}\B{1}$

12 $\Ees{\B^{}}{\W^{2}}{\B^{1}}{\B^{}}{\W^{3}}{\W^{2}}{\B^{}}\B{1}$

13 $\Ees{\B^{}}{\W^{0}}{\B^{}}{\B^{}}{\W^{3}}{\W^{2}}{\B^{}}\B{1}$

14 $\Ees{\B^{}}{\W^{2}}{\B^{}}{\B^{}}{\W^{5}}{\W^{3}}{\B^{1}}\B{1}$

15 $\Ees{\B^{}}{\W^{0}}{\B^{}}{\B^{1}}{\W^{2}}{\W^{2}}{\B^{}}\B{}$

}}{
\vtop{\hsize=0.25\hsize\parindent0pt

16 $\Ees{\B^{1}}{\W^{1}}{\B^{1}}{\B^{}}{\W^{2}}{\W^{2}}{\B^{}}\B{}$

17 $\Ees{\B^{1}}{\W^{1}}{\B^{}}{\B^{}}{\W^{4}}{\W^{3}}{\B^{1}}\B{}$

18 $\Ees{\B^{}}{\W^{2}}{\B^{}}{\B^{1}}{\W^{4}}{\W^{3}}{\B^{1}}\B{}$

20 $\Ees{\B^{}}{\W^{2}}{\B^{}}{\B^{}}{\W^{2}}{\W^{2}}{\B^{}}\B{}$

\bigskip\noindent

a $\Ees{\B^{}}{\W^{0}}{\B^{}}{\B^{1}}{\W^{2}}{\W^{1}}{\B^{1}}\B{}$

b $\Ees{\B^{}}{\W^{0}}{\B^{}}{\B^{}}{\W^{3}}{\W^{2}}{\B^{}}\B{1}$

c $\Ees{\B^{}}{\W^{0}}{\B^{}}{\B^{1}}{\W^{2}}{\W^{2}}{\B^{}}\B{}$

d $\Ees{\B^{}}{\W^{2}}{\B^{}}{\B^{}}{\W^{2}}{\W^{1}}{\B^{1}}\B{}$

e $\Ees{\B^{}}{\W^{2}}{\B^{1}}{\B^{}}{\W^{3}}{\W^{2}}{\B^{}}\B{1}$

}}{
\vtop{\hsize=0.25\hsize\parindent0pt\noindent

f $\Ees{\B^{1}}{\W^{1}}{\B^{1}}{\B^{}}{\W^{2}}{\W^{2}}{\B^{}}\B{}$

g $\Ees{\B^{}}{\W^{2}}{\B^{}}{\B^{1}}{\W^{4}}{\W^{2}}{\B^{}}\B{}$

h $\Ees{\B^{1}}{\W^{1}}{\B^{}}{\B^{}}{\W^{4}}{\W^{2}}{\B^{}}\B{}$

i $\Ees{\B^{}}{\W^{2}}{\B^{}}{\B^{}}{\W^{5}}{\W^{3}}{\B^{1}}\B{1}$

j $\Ees{\B^{1}}{\W^{1}}{\B^{}}{\B^{}}{\W^{4}}{\W^{3}}{\B^{1}}\B{}$

k $\Ees{\B^{}}{\W^{2}}{\B^{}}{\B^{1}}{\W^{4}}{\W^{3}}{\B^{1}}\B{}$

l $\Ees{\B^{}}{\W^{2}}{\B^{}}{\B^{}}{\W^{2}}{\W^{2}}{\B^{}}\B{}$

m $\Ees{\B^{}}{\W^{2}}{\B^{1}}{\B^{}}{\W^{6}}{\W^{4}}{\B^{}}\B{}$

n $\Ees{\B^{1}}{\W^{3}}{\B^{}}{\B^{}}{\W^{6}}{\W^{4}}{\B^{}}\B{}$

o $\Ees{\B^{}}{\W^{2}}{\B^{}}{\B^{}}{\W^{5}}{\W^{4}}{\B^{}}\B{1}$

}}}

\noindent
$A\ra{1}B$\quad
$a\ra{1}C\ra{1}d$\quad
$b\ra{1}D\ra{1}e$\quad
$c\ra{1}f$\quad
$f\ra{1}l$\quad
$h\ra{1}g$\quad
$j\ra{1}k$\quad
$m\ra{1}n$\newline
$A\ra{2}C$\quad
$B\ra{2}d$\quad
$b\ra{2}h$\quad
$c\ra{2}D\ra{2}g$\quad
$f\ra{2}e$\quad
$j\ra{2}i$\quad
$o\ra{2}n$\newline
$A\ra{3}a$\quad
$B\ra{3}C$\quad
$c\ra{3}b$\quad
$f\ra{3}D\ra{3}h$\quad
$k\ra{3}i$\quad
$l\ra{3}e\ra{3}g$\quad
$o\ra{3}m$\newline
$A\ra{4}D\ra{4}i$\quad
$B\ra{4}e$\quad
$C\ra{4}g$\quad
$a\ra{4}h$\quad
$c\ra{4}j\ra{4}m$\quad
$f\ra{4}k\ra{4}n$\newline
$A\ra{5}b$\quad
$B\ra{5}D$\quad
$C\ra{5}h$\quad
$d\ra{5}g$\quad
$e\ra{5}i$\quad
$f\ra{5}j$\quad
$l\ra{5}k\ra{5}m$\newline
$A\ra{6}c$\quad
$B\ra{6}f$\quad
$C\ra{6}D\ra{6}j$\quad
$a\ra{6}b$\quad
$d\ra{6}e\ra{6}k\ra{6}o$\quad
$g\ra{6}i$\quad
$i\ra{6}m$\newline
$A\ra{7}e$\quad
$a\ra{7}g$\quad
$b\ra{7}i$\quad
$c\ra{7}k$\quad
$j\ra{7}n$\newline
$A\ra{8}f$\quad
$B\ra{8}l$\quad
$C\ra{8}e$\quad
$a\ra{8}D\ra{8}k$\quad
$b\ra{8}j\ra{8}o$\quad
$h\ra{8}i\ra{8}n$\newline
$A\ra{9}i$\quad
$c\ra{9}m$\quad
$f\ra{9}n$\quad\quad
$B\ra{10}i$\quad
$f\ra{10}m$\quad
$l\ra{10}n$\newline
$A\ra{11}j$\quad
$B\ra{11}k$\quad
$C\ra{11}i$\quad
$D\ra{11}m$\quad
$e\ra{11}n$\quad
$f\ra{11}o$\newline
$A\ra{12}k$\quad
$D\ra{12}n$\quad
$a\ra{12}i$\quad
$b\ra{12}m$\quad
$c\ra{12}o$\quad\quad
$B\ra{13}j$\quad
$d\ra{13}i$\quad
$e\ra{13}m$\quad
$l\ra{13}o$\newline
$A\ra{14}m$\quad
$B\ra{14}n$\quad\quad
$C\ra{15}j$\quad
$d\ra{15}k$\quad
$e\ra{15}o$\quad
$g\ra{15}m$\newline
$C\ra{16}k$\quad
$D\ra{16}o$\quad
$a\ra{16}j$\quad
$g\ra{16}n$\quad
$h\ra{16}m$\quad\quad
$B\ra{17}o$\quad
$C\ra{17}m$\quad
$d\ra{17}n$\newline
$A\ra{18}o$\quad
$C\ra{18}n$\quad
$a\ra{18}m$\quad\quad
$a\ra{20}k$\quad
$b\ra{20}o$\quad
$h\ra{20}n$

Consider the coadjoint case, let $x\in\n$.
Then
$A\LRA{1}{4.B}B$,
$A\LRA{2}{1.B}C$,
$A\LRA{4}{5.A}D$.
Let $x_A$ be ample, then
$B\LRA{3}{4.B}C$,
$B\LRA{5}{4.B}D$.
Let $x_B$ be ample, 
then 
$C\LRA{6}{5.A}D$.

Adjoint case, $x\in\g/\p^-$.
Then 
$n\LRA{-12}{1.A}D$,
$n\LRA{-14}{1.A}B$,
$n\LRA{-18}{1.A}C$.
Let $x_n=0$, then
$m\LRA{-11}{4.A}D$,
$m\LRA{-14}{4.A}A$,
$m\LRA{-17}{4.A}C$.
Let $x_m=0$, then
$o\LRA{-16}{1.A}D$,
$o\LRA{-17}{1.A}B$,
$o\LRA{-18}{1.A}A$.
Let $x_o=0$, then
$i\LRA{-4}{1.A}D$,
$i\LRA{-11}{1.A}C$,
$i\LRA{-9}{4.B}A$.
Let $x_i=0$, then
$k\LRA{-8}{4.A}D$,
$k\LRA{-11}{1.A}B$,
$k\LRA{-12}{5.C}A$.
Let $x_k=0$, then
$j\LRA{-6}{1.A}D$,
$j\LRA{-11}{1.A}A$,
$j\LRA{-15}{4.B}C$.
Let $x_j=0$. If  $x_g\ne0$,
then by Proposition~5.E, we can apply $\exp(p_{-3})$ to force
$x_e$ to have rank $2$.
Afterwards, 
$e\LRA{-1}{5.C}D$,
$e\LRA{-4}{1.A}B$,
$e\LRA{-7}{1.C}A$.
Let $x_g=0$, then
$e\LRA{-1}{5.C}D$,
$e\LRA{-8}{4.A}C$,
$e\LRA{-4}{1.A}B$.
Let $x_e=0$, then
$h\LRA{-3}{1.A}D$,
$h\LRA{-5}{1.A}C$;
if $x_d\ne0$ then
$d\LRA{-2}{1.A}B$, if 
$x_d=0$ then $B\LRA{-1}{4.B}A$.
Let $x_h=0$, then
$D\LRA{-6}{5.A}C$,
$D\LRA{-5}{4.B}B$,
$D\LRA{-4}{5.A}A$.
Let $x_D$ be ample, then
$d\LRA{-1}{1.A}C$,
$d\LRA{-2}{1.A}B$.
Let $x_d=0$, then
$C\LRA{-3}{4.B}B$,
$C\LRA{-2}{1.B}A$.
Let $x_C$ be ample, 
then 
$B\LRA{-1}{4.B}A$.

\subhead Case 5D\endsubhead
In this case we study parabolic subgroup 16.

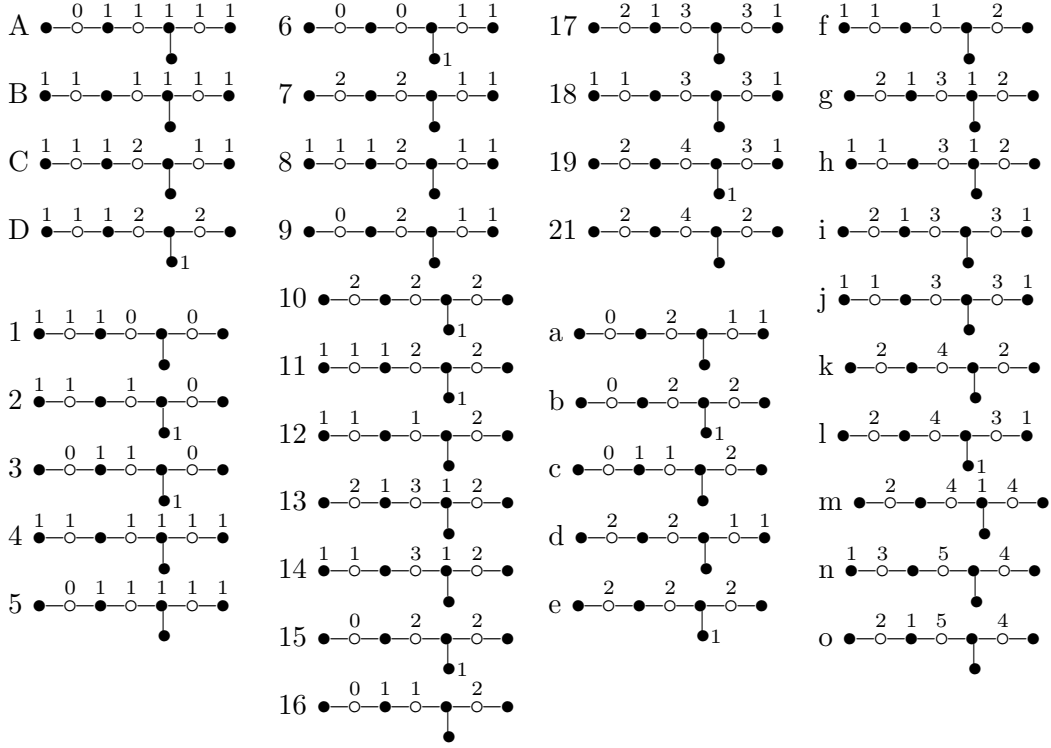
\hbox{{\vtop{\hsize=0.25\hsize\parindent0pt

A $\Ees{\B^{}}{\W^{0}}{\B^{1}}{\W^{1}}{\B^{1}}{\W^{1}}{\B^{1}}\B{}$

B $\Ees{\B^{1}}{\W^{1}}{\B^{}}{\W^{1}}{\B^{1}}{\W^{1}}{\B^{1}}\B{}$

C $\Ees{\B^{1}}{\W^{1}}{\B^{1}}{\W^{2}}{\B^{}}{\W^{1}}{\B^{1}}\B{}$

D $\Ees{\B^{1}}{\W^{1}}{\B^{1}}{\W^{2}}{\B^{}}{\W^{2}}{\B^{}}\B{1}$

\bigskip\noindent

1 $\Ees{\B^{1}}{\W^{1}}{\B^{1}}{\W^{0}}{\B^{}}{\W^{0}}{\B^{}}\B{}$

2 $\Ees{\B^{1}}{\W^{1}}{\B^{}}{\W^{1}}{\B^{}}{\W^{0}}{\B^{}}\B{1}$

3 $\Ees{\B^{}}{\W^{0}}{\B^{1}}{\W^{1}}{\B^{}}{\W^{0}}{\B^{}}\B{1}$

4 $\Ees{\B^{1}}{\W^{1}}{\B^{}}{\W^{1}}{\B^{1}}{\W^{1}}{\B^{1}}\B{}$

5 $\Ees{\B^{}}{\W^{0}}{\B^{1}}{\W^{1}}{\B^{1}}{\W^{1}}{\B^{1}}\B{}$

}}{
\vtop{\hsize=0.25\hsize\parindent0pt

6 $\Ees{\B^{}}{\W^{0}}{\B^{}}{\W^{0}}{\B^{}}{\W^{1}}{\B^{1}}\B{1}$

7 $\Ees{\B^{}}{\W^{2}}{\B^{}}{\W^{2}}{\B^{}}{\W^{1}}{\B^{1}}\B{}$

8 $\Ees{\B^{1}}{\W^{1}}{\B^{1}}{\W^{2}}{\B^{}}{\W^{1}}{\B^{1}}\B{}$

9 $\Ees{\B^{}}{\W^{0}}{\B^{}}{\W^{2}}{\B^{}}{\W^{1}}{\B^{1}}\B{}$

10 $\Ees{\B^{}}{\W^{2}}{\B^{}}{\W^{2}}{\B^{}}{\W^{2}}{\B^{}}\B{1}$

11 $\Ees{\B^{1}}{\W^{1}}{\B^{1}}{\W^{2}}{\B^{}}{\W^{2}}{\B^{}}\B{1}$

12 $\Ees{\B^{1}}{\W^{1}}{\B^{}}{\W^{1}}{\B^{}}{\W^{2}}{\B^{}}\B{}$

13 $\Ees{\B^{}}{\W^{2}}{\B^{1}}{\W^{3}}{\B^{1}}{\W^{2}}{\B^{}}\B{}$

14 $\Ees{\B^{1}}{\W^{1}}{\B^{}}{\W^{3}}{\B^{1}}{\W^{2}}{\B^{}}\B{}$

15 $\Ees{\B^{}}{\W^{0}}{\B^{}}{\W^{2}}{\B^{}}{\W^{2}}{\B^{}}\B{1}$

16 $\Ees{\B^{}}{\W^{0}}{\B^{1}}{\W^{1}}{\B^{}}{\W^{2}}{\B^{}}\B{}$

}}{
\vtop{\hsize=0.25\hsize\parindent0pt

17 $\Ees{\B^{}}{\W^{2}}{\B^{1}}{\W^{3}}{\B^{}}{\W^{3}}{\B^{1}}\B{}$

18 $\Ees{\B^{1}}{\W^{1}}{\B^{}}{\W^{3}}{\B^{}}{\W^{3}}{\B^{1}}\B{}$

19 $\Ees{\B^{}}{\W^{2}}{\B^{}}{\W^{4}}{\B^{}}{\W^{3}}{\B^{1}}\B{1}$

21 $\Ees{\B^{}}{\W^{2}}{\B^{}}{\W^{4}}{\B^{}}{\W^{2}}{\B^{}}\B{}$

\bigskip\noindent

a $\Ees{\B^{}}{\W^{0}}{\B^{}}{\W^{2}}{\B^{}}{\W^{1}}{\B^{1}}\B{}$

b $\Ees{\B^{}}{\W^{0}}{\B^{}}{\W^{2}}{\B^{}}{\W^{2}}{\B^{}}\B{1}$

c $\Ees{\B^{}}{\W^{0}}{\B^{1}}{\W^{1}}{\B^{}}{\W^{2}}{\B^{}}\B{}$

d $\Ees{\B^{}}{\W^{2}}{\B^{}}{\W^{2}}{\B^{}}{\W^{1}}{\B^{1}}\B{}$

e $\Ees{\B^{}}{\W^{2}}{\B^{}}{\W^{2}}{\B^{}}{\W^{2}}{\B^{}}\B{1}$

}}{
\vtop{\hsize=0.25\hsize\parindent0pt\noindent

f $\Ees{\B^{1}}{\W^{1}}{\B^{}}{\W^{1}}{\B^{}}{\W^{2}}{\B^{}}\B{}$

g $\Ees{\B^{}}{\W^{2}}{\B^{1}}{\W^{3}}{\B^{1}}{\W^{2}}{\B^{}}\B{}$

h $\Ees{\B^{1}}{\W^{1}}{\B^{}}{\W^{3}}{\B^{1}}{\W^{2}}{\B^{}}\B{}$

i $\Ees{\B^{}}{\W^{2}}{\B^{1}}{\W^{3}}{\B^{}}{\W^{3}}{\B^{1}}\B{}$

j $\Ees{\B^{1}}{\W^{1}}{\B^{}}{\W^{3}}{\B^{}}{\W^{3}}{\B^{1}}\B{}$

k $\Ees{\B^{}}{\W^{2}}{\B^{}}{\W^{4}}{\B^{}}{\W^{2}}{\B^{}}\B{}$

l $\Ees{\B^{}}{\W^{2}}{\B^{}}{\W^{4}}{\B^{}}{\W^{3}}{\B^{1}}\B{1}$

m $\Ees{\B^{}}{\W^{2}}{\B^{}}{\W^{4}}{\B^{1}}{\W^{4}}{\B^{}}\B{}$

n $\Ees{\B^{1}}{\W^{3}}{\B^{}}{\W^{5}}{\B^{}}{\W^{4}}{\B^{}}\B{}$

o $\Ees{\B^{}}{\W^{2}}{\B^{1}}{\W^{5}}{\B^{}}{\W^{4}}{\B^{}}\B{}$

}}}

\noindent
$A\ra{1}B$\quad
$a\ra{1}C\ra{1}d$\quad
$b\ra{1}D\ra{1}e$\quad
$c\ra{1}f$\quad
$h\ra{1}g$\quad
$j\ra{1}i$\quad
$o\ra{1}n$\newline
$A\ra{2}C$\quad
$B\ra{2}d$\quad
$b\ra{2}h\ra{2}k$\quad
$c\ra{2}D\ra{2}g$\quad
$f\ra{2}e$\quad
$j\ra{2}l$\quad
$m\ra{2}n$\newline
$A\ra{3}a$\quad
$B\ra{3}C$\quad
$c\ra{3}b$\quad
$e\ra{3}g\ra{3}k$\quad
$f\ra{3}D\ra{3}h$\quad
$i\ra{3}l$\quad
$m\ra{3}o$\newline
$A\ra{4}D\ra{4}i$\quad
$B\ra{4}e$\quad
$C\ra{4}g$\quad
$a\ra{4}h\ra{4}l\ra{4}n$\quad
$b\ra{4}j\ra{4}m$\newline
$A\ra{5}b$\quad
$B\ra{5}D\ra{5}j$\quad
$C\ra{5}h$\quad
$d\ra{5}g\ra{5}l\ra{5}o$\quad
$e\ra{5}i\ra{5}m$\newline
$A\ra{6}c$\quad
$B\ra{6}f$\quad
$C\ra{6}D$\quad
$a\ra{6}b$\quad
$d\ra{6}e$\quad
$g\ra{6}i$\quad
$h\ra{6}j$\quad
$k\ra{6}l\ra{6}m$\newline
$A\ra{7}g$\quad
$a\ra{7}k$\quad
$b\ra{7}l$\quad
$c\ra{7}i$\quad
$j\ra{7}n$\newline
$A\ra{8}h$\quad
$B\ra{8}g$\quad
$C\ra{8}k$\quad
$D\ra{8}l$\quad
$c\ra{8}j\ra{8}o$\quad
$f\ra{8}i\ra{8}n$\newline
$B\ra{9}h$\quad
$d\ra{9}k$\quad
$e\ra{9}l$\quad
$f\ra{9}j$\quad
$i\ra{9}o$\newline
$A\ra{10}i$\quad
$a\ra{10}l$\quad
$b\ra{10}m$\quad
$h\ra{10}n$\newline
$A\ra{11}j$\quad
$B\ra{11}i$\quad
$C\ra{11}l$\quad
$D\ra{11}m$\quad
$g\ra{11}n$\quad
$h\ra{11}o$\newline
$C\ra{12}i$\quad
$a\ra{12}j$\quad
$h\ra{12}m$\quad
$k\ra{12}n$\quad\quad
$A\ra{13}l$\quad
$D\ra{13}n$\quad
$b\ra{13}o$\quad
$c\ra{13}m$\newline
$B\ra{14}l$\quad
$D\ra{14}o$\quad
$e\ra{14}n$\quad
$f\ra{14}m$\quad\quad
$B\ra{15}j$\quad
$d\ra{15}l$\quad
$e\ra{15}m$\quad
$g\ra{15}o$\newline
$C\ra{16}j$\quad
$d\ra{16}i$\quad
$g\ra{16}m$\quad
$k\ra{16}o$\quad\quad
$A\ra{17}m$\quad
$C\ra{17}n$\quad
$a\ra{17}o$\newline
$B\ra{18}m$\quad
$C\ra{18}o$\quad
$d\ra{18}n$\quad\quad
$A\ra{19}o$\quad
$B\ra{19}n$\quad\quad
$c\ra{21}o$\quad
$f\ra{21}n$

Coadjoint case, $x\in\n$.
Then $A\LRA{1}{1.B}B$,
$A\LRA{2}{4.B}C$, $A\LRA{4}{5.A}D$.
Let $x_A$ be ample. Then
$B\LRA{3}{4.B}C$, 
$B\LRA{5}{5.A}D$.
Let $x_B$ be ample. Then
$C\LRA{6}{4.B}D$.

Adjoint case, $x\in\g/\p^-$ 
Then 
$n\LRA{-13}{1.A}D$,
$n\LRA{-17}{1.A}C$,
$n\LRA{-19}{1.A}B$.
Let $x_n=0$. Then
$o\LRA{-14}{1.A}D$,
$o\LRA{-18}{1.A}C$,
$o\LRA{-19}{1.A}A$.
Let $x_o=0$. Then
$m\LRA{-11}{4.A}D$,
$m\LRA{-17}{4.A}A$,
$m\LRA{-18}{4.A}B$.
Let $x_m=0$. Then
$l\LRA{-8}{4.A}D$,
$l\LRA{-11}{1.A}C$,
$l\LRA{-13}{5.C}A$.
Let $x_l=0$. Then
$i\LRA{-4}{1.A}D$,
$i\LRA{-11}{1.A}B$,
$i\LRA{-10}{4.B}A$.
Let $x_i=0$. Then
$j\LRA{-5}{1.A}D$,
$j\LRA{-11}{1.A}A$,
$j\LRA{-15}{4.B}B$.
Let $x_j=0$, $x_k\ne0$.
We can apply $\exp(p_{-3})$ to force
$x_g$ to have rank $2$,
then 
$g\LRA{-2}{5.C}D$,
$g\LRA{-4}{1.A}C$,
$g\LRA{-7}{1.C}A$.
Let $x_k=0$, then
$g\LRA{-2}{5.C}D$,
$g\LRA{-4}{1.A}C$,
$g\LRA{-8}{4.A}B$.
Let $x_g=0$, then
$h\LRA{-3}{5.C}D$,
$h\LRA{-5}{1.A}C$,
$h\LRA{-8}{4.A}A$.
Let $x_h=0$, $x_e\ne0$.
Applying $\exp(p_{-6})$ if necessary, we make $x_d=0$,
then 
$e\LRA{-1}{4.A}D$,
$e\LRA{-4}{4.A}B$,
$C\LRA{-2}{4.B}A$.
Let $x_e=0$. Then
$D\LRA{-6}{4.B}C$,
$D\LRA{-5}{5.A}B$,
$D\LRA{-4}{5.A}A$.
Suppose $x_D$ is ample.
Then
$d\LRA{-1}{1.A}C$,
$d\LRA{-2}{1.A}B$.
Let $x_d=0$. Then
$C\LRA{-3}{4.B}B$,
$C\LRA{-2}{4.B}A$.
Let  $x_C$ be ample, then
$B\LRA{-1}{1.B}A$.

\subhead Case 5E\endsubhead
This case consists of parabolic subgroups 17 and 18.
We shall consider number 17.

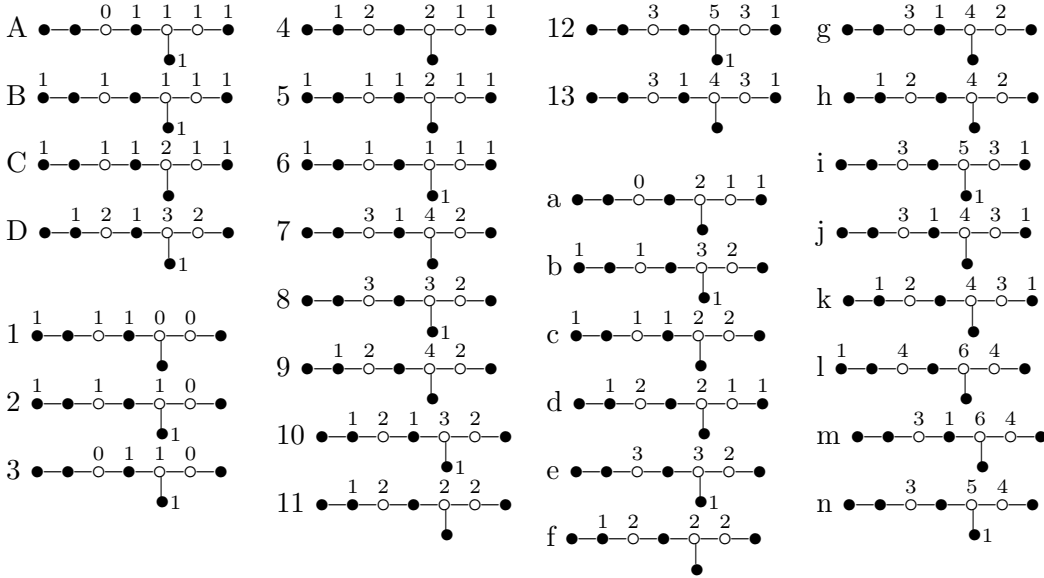
\hbox{{\vtop{\hsize=0.25\hsize\parindent0pt

A $\Ees{\B^{}}{\B^{}}{\W^{0}}{\B^{1}}{\W^{1}}{\W^{1}}{\B^{1}}\B{1}$

B $\Ees{\B^{1}}{\B^{}}{\W^{1}}{\B^{}}{\W^{1}}{\W^{1}}{\B^{1}}\B{1}$

C $\Ees{\B^{1}}{\B^{}}{\W^{1}}{\B^{1}}{\W^{2}}{\W^{1}}{\B^{1}}\B{}$

D $\Ees{\B^{}}{\B^{1}}{\W^{2}}{\B^{1}}{\W^{3}}{\W^{2}}{\B^{}}\B{1}$

\bigskip\noindent

1 $\Ees{\B^{1}}{\B^{}}{\W^{1}}{\B^{1}}{\W^{0}}{\W^{0}}{\B^{}}\B{}$

2 $\Ees{\B^{1}}{\B^{}}{\W^{1}}{\B^{}}{\W^{1}}{\W^{0}}{\B^{}}\B{1}$

3 $\Ees{\B^{}}{\B^{}}{\W^{0}}{\B^{1}}{\W^{1}}{\W^{0}}{\B^{}}\B{1}$

}}{
\vtop{\hsize=0.25\hsize\parindent0pt

4 $\Ees{\B^{}}{\B^{1}}{\W^{2}}{\B^{}}{\W^{2}}{\W^{1}}{\B^{1}}\B{}$

5 $\Ees{\B^{1}}{\B^{}}{\W^{1}}{\B^{1}}{\W^{2}}{\W^{1}}{\B^{1}}\B{}$

6 $\Ees{\B^{1}}{\B^{}}{\W^{1}}{\B^{}}{\W^{1}}{\W^{1}}{\B^{1}}\B{1}$

7 $\Ees{\B^{}}{\B^{}}{\W^{3}}{\B^{1}}{\W^{4}}{\W^{2}}{\B^{}}\B{}$

8 $\Ees{\B^{}}{\B^{}}{\W^{3}}{\B^{}}{\W^{3}}{\W^{2}}{\B^{}}\B{1}$

9 $\Ees{\B^{}}{\B^{1}}{\W^{2}}{\B^{}}{\W^{4}}{\W^{2}}{\B^{}}\B{}$

10 $\Ees{\B^{}}{\B^{1}}{\W^{2}}{\B^{1}}{\W^{3}}{\W^{2}}{\B^{}}\B{1}$

11 $\Ees{\B^{}}{\B^{1}}{\W^{2}}{\B^{}}{\W^{2}}{\W^{2}}{\B^{}}\B{}$

}}{
\vtop{\hsize=0.25\hsize\parindent0pt

12 $\Ees{\B^{}}{\B^{}}{\W^{3}}{\B^{}}{\W^{5}}{\W^{3}}{\B^{1}}\B{1}$

13 $\Ees{\B^{}}{\B^{}}{\W^{3}}{\B^{1}}{\W^{4}}{\W^{3}}{\B^{1}}\B{}$

\bigskip\noindent

a $\Ees{\B^{}}{\B^{}}{\W^{0}}{\B^{}}{\W^{2}}{\W^{1}}{\B^{1}}\B{}$

b $\Ees{\B^{1}}{\B^{}}{\W^{1}}{\B^{}}{\W^{3}}{\W^{2}}{\B^{}}\B{1}$

c $\Ees{\B^{1}}{\B^{}}{\W^{1}}{\B^{1}}{\W^{2}}{\W^{2}}{\B^{}}\B{}$

d $\Ees{\B^{}}{\B^{1}}{\W^{2}}{\B^{}}{\W^{2}}{\W^{1}}{\B^{1}}\B{}$

e $\Ees{\B^{}}{\B^{}}{\W^{3}}{\B^{}}{\W^{3}}{\W^{2}}{\B^{}}\B{1}$

f $\Ees{\B^{}}{\B^{1}}{\W^{2}}{\B^{}}{\W^{2}}{\W^{2}}{\B^{}}\B{}$

}}{
\vtop{\hsize=0.25\hsize\parindent0pt\noindent

g $\Ees{\B^{}}{\B^{}}{\W^{3}}{\B^{1}}{\W^{4}}{\W^{2}}{\B^{}}\B{}$

h $\Ees{\B^{}}{\B^{1}}{\W^{2}}{\B^{}}{\W^{4}}{\W^{2}}{\B^{}}\B{}$

i $\Ees{\B^{}}{\B^{}}{\W^{3}}{\B^{}}{\W^{5}}{\W^{3}}{\B^{1}}\B{1}$

j $\Ees{\B^{}}{\B^{}}{\W^{3}}{\B^{1}}{\W^{4}}{\W^{3}}{\B^{1}}\B{}$

k $\Ees{\B^{}}{\B^{1}}{\W^{2}}{\B^{}}{\W^{4}}{\W^{3}}{\B^{1}}\B{}$

l $\Ees{\B^{1}}{\B^{}}{\W^{4}}{\B^{}}{\W^{6}}{\W^{4}}{\B^{}}\B{}$

m $\Ees{\B^{}}{\B^{}}{\W^{3}}{\B^{1}}{\W^{6}}{\W^{4}}{\B^{}}\B{}$

n $\Ees{\B^{}}{\B^{}}{\W^{3}}{\B^{}}{\W^{5}}{\W^{4}}{\B^{}}\B{1}$

}}}

\noindent 
$A\ra{1}B$\quad
$a\ra{1}C\ra{1}d$\quad
$b\ra{1}D\ra{1}e$\quad
$c\ra{1}f$\quad
$h\ra{1}g$\quad
$k\ra{1}j$\quad
$m\ra{1}l$\newline
$A\ra{2}C$\quad
$B\ra{2}d$\quad
$b\ra{2}h$\quad
$c\ra{2}D\ra{2}g$\quad
$f\ra{2}e$\quad
$k\ra{2}i$\quad
$n\ra{2}l$\newline
$A\ra{3}a$\quad
$B\ra{3}C$\quad
$c\ra{3}b$\quad
$e\ra{3}g$\quad
$f\ra{3}D\ra{3}h$\quad
$j\ra{3}i$\quad
$n\ra{3}m$\newline
$A\ra{4}D$\quad
$B\ra{4}e$\quad
$C\ra{4}g$\quad
$a\ra{4}h$\quad
$b\ra{4}i$\quad
$c\ra{4}j$\quad
$k\ra{4}l$\newline
$A\ra{5}b$\quad
$B\ra{5}D\ra{5}i$\quad
$C\ra{5}h$\quad
$c\ra{5}k\ra{5}m$\quad
$d\ra{5}g$\quad
$f\ra{5}j\ra{5}l$\newline
$A\ra{6}c$\quad
$B\ra{6}f$\quad
$C\ra{6}D\ra{6}j$\quad
$a\ra{6}b\ra{6}k\ra{6}n$\quad
$d\ra{6}e$\quad
$h\ra{6}i\ra{6}l$\newline
$A\ra{7}i$\quad
$c\ra{7}l$\quad\quad
$A\ra{8}j$\quad
$a\ra{8}i$\quad
$b\ra{8}l$\newline
$B\ra{9}i$\quad
$c\ra{9}m$\quad
$f\ra{9}l$\quad\quad
$A\ra{10}k$\quad
$B\ra{10}j$\quad
$C\ra{10}i$\quad
$D\ra{10}l$\quad
$b\ra{10}m$\quad
$c\ra{10}n$\newline
$C\ra{11}j$\quad
$a\ra{11}k$\quad
$b\ra{11}n$\quad
$h\ra{11}l$\quad\quad
$A\ra{12}m$\quad
$B\ra{12}l$\newline
$A\ra{13}n$\quad
$C\ra{13}l$\quad
$a\ra{13}m$

Consider the coadjoint case, let $x\in\n$.
Then $A\LRA{1}{4.B}B$, $A\LRA{2}{4.B}C$, $A\LRA{4}{4.B}D$.
Let $x_A$ be ample.
Then $B\LRA{3}{1.B}C$, $B\LRA{5}{5.A}D$.
Suppose that $x_B$ is also ample.
Then $C\LRA{6}{5.A}D$.

Adjoint case, let $x\in\g/\p^-$ 
Then $l\LRA{-10}{4.A}D$, $l\LRA{-12}{4.A}B$,
$l\LRA{-13}{4.A}C$.
Let $x_l=0$. Then $i\LRA{-5}{1.A}D$,
$i\LRA{-10}{1.A}C$, $i\LRA{-9}{4.B}B$.
Let $x_i=0$. Then $j\LRA{-6}{1.A}D$,
$j\LRA{-10}{1.A}B$,
$j\LRA{-11}{4.B}C$..
Let $x_j=0$, $x_g\ne0$.
Then we first apply 
$\exp(p_{-3})$ and $\exp(p_{-5})$ to make
$x_e=x_d=0$, then $g\LRA{-2}{1.A}D$, $g\LRA{-4}{1.A}C$,
$B\LRA{-1}{4.B}A$.
Let $x_g=0$, $x_e\ne0$.
First, we apply $\exp(p_{-6})$ to make $x_d=0$, then
$e\LRA{-1}{1.A}D$, $e\LRA{-4}{1.A}B$, 
$C\LRA{-2}{4.B}A$.
Let $x_e=0$. Then $d\LRA{-1}{5.C}C$,
$d\LRA{-2}{5.C}B$, $D\LRA{-4}{4.B}A$.
Let $x_d=0$. Then $h\LRA{-3}{4.A}D$,
$h\LRA{-5}{4.A}C$, $B\LRA{-1}{4.B}A$.
Let $x_h=0$. Then $D\LRA{-6}{5.A}C$,
$D\LRA{-5}{5.A}B$, $D\LRA{-4}{4.B}A$.
Let $x_D$ be ample. Then $C\LRA{-3}{1.B}B$,
$C\LRA{-2}{4.B}A$.
Let $x_C$ be ample, then 
$B\LRA{-1}{4.B}A$.

\head \S5. Table\endhead

\centerline{\vbox{\offinterlineskip
\halign
{\vrule#&
  \strut\ \hfil# \hfil&        \vrule#&
  \strut\ \hfil#\ \qquad \hfil&        \vrule#&
  \strut\ \hfil# \hfil&        \vrule#&
  \strut\ \hfil#\ \qquad \hfil&        \vrule#&
  \strut\ \hfil# \hfil&        \vrule#&
  \strut\ \hfil#\ \qquad \hfil&        \vrule#&
  \strut\ \hfil# \hfil&        \vrule#&
  \strut\ \hfil#\ \qquad \hfil&        \vrule#
 \cr
\noalign{\hrule}\cr
& 1 && $\Es\B\W\B\W\B\W\B$ &
& 2 && $\Es\B\W\B\W\W\B\B$ &
& 3 && $\Ee\W\B\W\B\W\W\B\B$ &
& 4 && $\Ee\W\B\W\B\W\B\W\B$ &\cr
\noalign{\hrule}\cr
& 5 && $\Ee\B\W\W\B\W\W\B\B$ &
& 6 && $\Ee\B\W\W\B\W\B\W\B$ &
& 7 && $\Ee\B\W\B\W\W\W\B\B$ &
& 8 && $\Ee\B\W\B\W\W\B\W\B$ &\cr
\noalign{\hrule}\cr
& 9 && $\Ee\B\W\B\W\B\W\B\W$ &
& 10 && $\Es\B\W\B\B\B\W\B$ &
& 11 && $\Ee\W\B\W\B\W\B\B\B$ &
& 12 && $\Ee\B\W\W\B\W\B\B\B$ &\cr
\noalign{\hrule}\cr
& 13 && $\Ee\B\W\B\W\W\B\B\B$ &
& 14 && $\Ee\B\W\B\B\W\W\B\B$ &
& 15 && $\Ee\B\W\B\B\W\B\W\B$ &
& 16 && $\Ee\B\W\B\W\B\W\B\B$ &\cr
\noalign{\hrule}\cr
& 17 && $\Ee\B\B\W\B\W\W\B\B$ &
& 18 && $\Ee\B\B\W\B\W\B\W\B$ &
& 19 && $\Es\B\W\B\W\B\B\B$ &
& 20 && $\Es\B\B\B\W\W\B\B$ &\cr
\noalign{\hrule}\cr
& 21 && $\Es\B\W\B\B\W\B\B$ &
& 22 && $\Ee\W\B\B\B\W\W\B\B$ &
& 23 && $\Ee\W\B\B\B\W\B\W\B$ &
& 24 && $\Ee\B\B\B\W\W\W\B\B$ &\cr
\noalign{\hrule}\cr
& 25 && $\Ee\B\B\B\W\W\B\W\B$ &
& 26 && $\Ee\B\B\B\W\B\W\B\W$ &
& 27 && $\Ee\W\B\W\B\B\W\B\B$ &
& 28 && $\Ee\B\W\W\B\B\W\B\B$ &\cr
\noalign{\hrule}\cr
& 29 && $\Ee\B\W\B\B\B\W\B\W$ &
& 30 && $\Ee\B\W\B\W\B\B\W\B$ &
& 31 && $\Ee\B\W\B\W\B\B\B\W$ &
& 32 && $\Ee\B\B\B\W\B\B\W\B$ &\cr
\noalign{\hrule}\cr
& 33 && $\Ee\B\B\B\W\B\B\B\W$ &
& 34 && $\Es\B\B\B\W\B\W\B$ &
& 35 && $\Ee\W\B\W\B\B\B\W\B$ &
& 36 && $\Ee\B\W\W\B\B\B\W\B$ &\cr
\noalign{\hrule}\cr
& 37 && $\Ee\B\B\W\B\B\B\W\B$ &
& 38 && $\Ee\B\W\B\B\W\B\B\B$ &
& 39 && $\Ee\B\B\W\B\W\B\B\B$ &
& 40 && $\Es\B\B\B\W\B\B\B$ &\cr
\noalign{\hrule}\cr
& 41 && $\Ee\B\B\W\B\B\W\B\B$ &
& 42 && $\Ee\B\B\B\W\B\W\B\B$ &
& 43 && $\Ee\W\B\B\B\W\B\B\B$ &
& 44 && $\Ee\B\B\B\W\W\B\B\B$ &\cr
\noalign{\hrule}\cr
& 45 && $\Ee\W\B\W\B\B\B\B\B$ &
& 46 && $\Ee\B\W\W\B\B\B\B\B$ &
& 47 && $\Ee\W\B\B\B\B\W\B\B$ &
& 48 && $\Ee\B\B\B\B\B\W\B\W$ &\cr
\noalign{\hrule}\cr
& 49 && $\Ee\W\B\B\B\B\B\W\B$ &
& 50 && $\Ee\B\W\B\W\B\B\B\B$ &
& 51 && $\Ee\B\W\B\B\B\W\B\B$ &
& 52 && $\Ee\B\W\B\B\B\B\W\B$ &\cr
\noalign{\hrule}\cr
& 53 && $\Ee\B\W\B\B\B\B\B\W$ &
& 54 && $\Ee\B\B\W\B\B\B\B\B$ &
& 55 && $\Ee\B\B\B\W\B\B\B\B$ &
& 56 && $\Ee\B\B\B\B\W\W\B\B$ &\cr
\noalign{\hrule}\cr
& 57 && $\Ee\B\B\B\B\W\B\W\B$ &
& 58 && $\Ee\B\B\B\B\W\B\B\B$ &
& 59 && $\Ee\B\B\B\B\B\W\B\B$ &&&&&\cr
\noalign{\hrule}\cr
}
}}
\bigskip

\head References\endhead
\widestnumber\key{XXXX}

\ref\key BHR
\by T. Br\"ustle, L. Hille, G. R\"ohrle
\paper Finiteness for parabolic group actions in classical groups
\jour Arch. Math.
\vol 76(2)
\yr 2001
\pages 81--87
\endref

\ref\key Ka
\by V.~G.~Kac
\paper Some remarks on nilpotent orbits
\jour J. Algebra
\vol 64
\yr 1980
\pages 190-213
\endref

\ref\key Li
\by W.~Lichtenstein
\paper A system of quadrics describing the orbit of the highest weight vector
\jour Proc.~ A.M.S.
\vol 84
\yr 1982
\pages 605--608
\endref

\ref\key LW
\by R. Lipsman, J. Wolf
\paper Canonical semi-invariants and the Plancherel formula for parabolic groups
\jour Trans. Amer. Math. Soc. 
\vol 269(1) 
\yr 1982
\pages 111--131
\endref

\ref\key PR
\by V. Popov, G. R\"ohrle
\paper On the number of orbits of a parabolic subgroup
on its unipotent radical
\jour Algebraic Groups and Lie Groups, Lehrer G.I., Ed., Austr. Math. Soc. Lecture 
Series
\vol 9
\yr 1997
\endref

\ref\key Ri
\by R.W.~Richardson
\paper Conjugacy classes of $n$-tuples in Lie algebras and
algebraic groups
\jour Duke Math. J. {\bf 57} (1988), 1--35
\endref

\ref\key Te1
\by E. Tevelev
\paper Moore--Penrose inverse, parabolic subgroups, and Jordan pairs
\jour to appear in:  Journal of Lie Theory
\endref

\ref\key Te2
\by E. Tevelev
\paper Isotropic subspaces of multi--linear forms
\jour Matematicheskie Zametki
\vol 69(6)
\yr 2001
\pages 925--933 (in Russian)
\endref

\ref\key VO
\by E.~Vinberg, A.~Onischik
\publ ``Seminar on Lie Groups and Algebraic Groups'',
Berlin: Springer, 1990
\endref

\ref\key Vi
\by E.~Vinberg
\paper On stability of actions of reductive algebraic groups. 
\jour In: Lie Algebras, Rings and Related Topics 
(Fong Yuen, A.A.Mikhalev, E.Zelmanov, Eds.),
Springer, Hong Kong 
\yr 2000
\pages 188-202
\endref

\enddocument